\documentclass[12pt]{article}
\usepackage{amsmath}
\usepackage{amssymb}
\usepackage{amsthm}
\usepackage{fullpage}
\usepackage{mathtools}
\usepackage{color}
\usepackage{tikz}
\usepackage{pgfplots}
  \usepackage{grffile}
\usepackage{combelow}
\usepackage{comment}
\usepackage{amsfonts}
\usepackage{color,ulem}
\usepackage{fullpage}
\usepackage{graphicx}%
\usepackage[active]{srcltx}
\usepackage{mathabx}
\usepackage{verbatim}
\usepackage{cancel}
\usepackage[official]{eurosym}

\usepackage{yhmath}

  \pgfplotsset{compat=newest}
  \usetikzlibrary{plotmarks}
  \usetikzlibrary{arrows.meta}
  \usepgfplotslibrary{patchplots}

\providecommand{\U}[1]{\protect\rule{.1in}{.1in}}

\newtheorem{theorem}{Theorem}

\newtheorem{lemma}[theorem]{Lemma}

\newtheorem{proposition}[theorem]{Proposition}
\newtheorem{remark}[theorem]{Remark}

\newcommand{\eps}{\epsilon}

\newcommand{\R}{\mathbb{R}}
\newcommand{\E}{\mathbb{E}}

\newcommand{\EE}{\mathbb{E}}

\newcommand{\C}{\mathcal{C}}

\newcommand{\dd}{\textup{d}}

\makeatletter
\newlength{\temp@wip@width}
\newlength{\temp@wip@height}
\newcommand{\wideinvparen}[1]{%
  \vfuzz=30pt% BAD: to remove overfull vbox warnings...
  \setlength{\temp@wip@width}{\widthof{$#1$}}%
  \setlength{\temp@wip@height}{\heightof{$#1$}}%
  #1\hspace{-\temp@wip@width}%
  \raisebox{\temp@wip@height+1pt}[\heightof{$\wideparen{#1}$}]%
    {\rotatebox[origin=c]{180}{\vbox to 0pt{\hbox{$\wideparen{\hphantom{#1}}$}}}}%
}

\def\parab{\wideparen}
\def\var{{\rm var}}
\def\cov{{\rm cov}}

\newcommand{\bh}{\boldsymbol{h}}
\newcommand{\bM}{\boldsymbol{M}}

\usepackage{accents}
\newlength{\dhatheight}

\usepackage[breaklinks,colorlinks=true,linkcolor=blue,citecolor=red,backref=page]{hyperref}

\begin{document}
\title{A control variate method  based on polynomial approximation of Brownian path}
\author{Josselin Garnier and Laurent Mertz}
\maketitle
\begin{abstract}
We present a novel control variate technique for enhancing the efficiency of Monte Carlo (MC) estimation of expectations involving solutions to stochastic differential equations (SDEs). Our method integrates a primary fine-time-step discretization of the SDE with a control variate derived from a secondary coarse-time-step discretization driven by a piecewise parabolic approximation of Brownian motion. This approximation is conditioned on the same fine-scale Brownian increments, enabling strong coupling between the estimators. The expectation of the control variate is computed via an independent MC simulation using the coarse approximation.
We characterize the minimized quadratic error decay as a function of the computational budget and the weak and strong orders of the primary and secondary 
discretization schemes.
We  demonstrate the method's effectiveness through numerical experiments on representative SDEs. %The proposed approach offers a scalable and accurate alternative for high-dimensional stochastic simulations.
%OLD
\begin{comment}
We propose a control variate method for efficient Monte Carlo (MC) estimation of the expectation of a function related to the solution of a stochastic differential equation (SDE). In a standard MC method for an SDE driven by Brownian motion, one generally uses a discretization method with fine time step to simulate the SDE multiple times and calculate an average outcome from these simulations. Here, we introduce a control variate coupled with the standard discretization method. The control variate is produced by an other discretized version of the same SDE but driven by a piecewise parabolic approximation of Brownian motion with a coarse time step. 
Crucially, this approximation of Brownian Motion is conditioned on the same small Brownian increments used in the discretization method with fine time step. The expectation of the control variate is estimated by an intensive MC sampling using the parabolic approximation of the Brownian motion without
conditioning, which is computationally much cheaper. 
The method's efficiency is demonstrated with several examples.
\end{comment}
\end{abstract}

%\tableofcontents

\section{Introduction}

The Monte Carlo (MC) method is one of the most effective numerical techniques for estimating the expectation associated with the solution of a stochastic differential equation (SDE) driven by Brownian motion, particularly in disciplines such as quantitative finance, physics, and engineering \cite{kloeden92,milstein04}. 
%In certain specific cases, alternative approaches to the MC method may rely on exact analytical formulas or partial differential equations (PDEs), notably via the Feynman-Kac formula. 
To obtain an accurate estimation, the method requires simulating many trajectories of a discretized version of the SDE and averaging the quantity of interest across these paths. Naturally, accuracy improves with a finer time discretization and a larger number of MC samples. 
However, to achieve high precision, the computational cost of such simulations can be substantial. 
To mitigate the cost, employing  high-order schemes and variance reduction techniques can be highly effective.

%variance reduction techniques --such as the control variate method \cite{glasserman} or the more recent multilevel Monte Carlo method \cite{giles08a}-- can be highly effective in this context.
%Using high-order schemes can also reduce the numerical cost and many schemes have been proposed \cite{clark80,rumelin82,talay90,newton91,newton94,milstein09,rossler10}.% including adaptive ones \cite{foster23}.

High-order discretization methods are often necessary to balance accuracy and efficiency. Standard schemes such as the Euler–Maruyama method converge with strong order $1/2$ and weak order $1$, which can be prohibitively inefficient.
%when high precision is required, especially in multi-dimensional problems or when computing sensitivities. 
Higher-order schemes, such as the Milstein method or stochastic Runge–Kutta (RK) methods, achieve better convergence rates by incorporating additional information from stochastic integrals \cite{kloeden92,milstein95}. This improved accuracy reduces the discretization bias, meaning fewer time steps are required to reach a desired tolerance, which is crucial when simulations are embedded within MC frameworks where computational costs scale with both the numbers of samples and time steps.

Variance reduction techniques have been extensively studied to overcome the slow convergence of classical MC methods. Early approaches focused on control variates (CVs) and antithetic variates \cite{glasserman}, as well as importance sampling strategies tailored to rare-event problems \cite{asmussen07}. 
CVs improve efficiency by exploiting correlations between the target functional and auxiliary quantities with known or cheaply estimable expectations, thereby reducing variance without biasing the estimator \cite{glasserman}. In SDE applications, natural CVs include analytically tractable approximations \cite{garnier22}, coarse discretization approximations \cite{giles08a}, and conditional probabilistic representations \cite{milstein09}. Recent developments extend this principle to adaptive and data-driven CVs, including regression-based methods \cite{milevsky98} and machine-learning–enhanced strategies \cite{hinds23}.
%In high-dimensional applications, quasi-Monte Carlo methods have been explored to exploit low-discrepancy sequences \cite{caflish88}, though their efficiency strongly depends on problem structure. 
A major breakthrough came with the introduction of the multilevel Monte Carlo (MLMC) method \cite{giles08a}, which exploits a hierarchy of discretizations to dramatically reduce variance at near-optimal computational complexity, spurring further extensions such as the multilevel quasi-MC \cite{giles09} and the multi-index MC framework \cite{haji16}.
%More recently, variance reduction has also been achieved by combining Monte Carlo with machine learning, such as the use of neural networks to construct problem-specific control variates \cite{bouchard22}. 

%The problem can be addressed from the opposite direction by starting with a computational budget constraint and determining how to achieve the highest possible accuracy within that budget.

In this paper, we introduce a novel CV approach designed under a computational budget constraint, which significantly enhances the efficiency of MC methods for SDEs driven by Brownian motion. Our method introduces a CV constructed from a discretized version of the same SDE, but driven by a piecewise parabolic approximation of Brownian motion \cite{lyons20}, using a coarse time step. The key innovation lies in conditioning this approximation on the same fine Brownian increments used in the primary SDE discretization with a fine time step. This coupling between the CV and the target quantity enables a significant reduction in variance. The expectation of the CV is estimated through MC sampling using the unconditioned parabolic approximation of BM, which is computationally much cheaper than simulating the full SDE with fine discretization. In this context, for a given computational budget, the four parameters -- fine and coarse time steps, and the two sizes of MC samples -- are optimized to achieve the highest accuracy.

The paper is organized as follows. Section~\ref{sec:overview} presents an overview of the results. It highlights the advantages of our method by comparing its minimized quadratic error decay rate with respect to the computational budget with that of the optimized standard MC method. Section~\ref{sec:main} presents our novel CV method, including the construction of both primary and secondary discretization schemes, the coupling mechanism, and the optimization of parameters under a fixed computational budget to achieve the maximal variance reduction. Section~\ref{sec:num} presents numerical illustrations that confirm the variance reduction, using several examples involving Euler and RK schemes for the primary discretization and the parabola method for the secondary scheme. Special attention is given to the inhomogeneous geometric Brownian motion (IGBM) \cite{tubikanec,zhao09}, which challenges the analysis due to the very small constant in the strong convergence error of the parabola method. In Section~\ref{sec:multi}, we discuss the integration of our approach into a multilevel MC framework.
%and outline directions for future research.

\section{Overview of the results}
\label{sec:overview}%
In our CV method, the central innovation lies in the construction of a CV using a piecewise parabolic approximation of Brownian motion \cite{lyons20}, conditioned on the fine-scale Brownian increments employed in the costly, primary SDE simulation. This CV is computed using a coarse time-step discretization, which is significantly less expensive than the fine-scale simulation. The method is formulated under a computational budget constraint, and we derive optimal choices for four key parameters: the fine and coarse time steps, and the respective numbers of MC samples.

Our theoretical analysis shows that the proposed CV estimator achieves a quadratic error decay rate of the form
\[
\mathcal{E} \lesssim \mathcal{C}^{-\frac{4\alpha\gamma + 2\alpha}{4\alpha\gamma + 2\alpha + 1}},
\]
%where $\alpha$ and $\gamma$ denote the weak and strong convergence orders of the discretization schemes.
where $\mathcal{C}$ is the computational budget, $\alpha$ denotes the weak  convergence order of the primary discretization scheme and $\gamma$ denotes the strong convergence error of the secondary scheme.
For instance, using the Euler--Maruyama method as the primary scheme ($\alpha = 1$) and the parabolic approximation of Brownian motion method as the secondary scheme ($\gamma = 1$) yields an error rate of $\mathcal{C}^{-6/7}$. This improves the optimized standard MC method based on the primary scheme where
\[
\mathcal{E} \lesssim \mathcal{C}^{-\frac{2\alpha}{2\alpha + 1}}.
\]
In the latter case, optimization is performed over only two parameters: the time step and the corresponding number of MC samples
\cite{duffie95}. For example, using the Euler--Maruyama method yields an error rate of $\mathcal{C}^{-2/3}$.
Schemes with higher order than the Euler--Maruyama method can be used, and similar improvements due to the use of a CV based on the parabolic approximation of Brownian motion
can be obtained.

Numerical experiments support the theoretical predictions and demonstrate the efficiency of the proposed CV method across a variety of SDE examples, including systems with multiplicative noise and Brownian particles in multi-well potentials. 
In these academic and realistic test cases, we will see that the theoretical predictions are confirmed for moderate and large computational budgets.
Particular attention is given to the inhomogeneous geometric Brownian motion, which is considered as a reference test case in the literature because it is a simple SDE but yet has no known method of exact simulation. We will see that this example involves constants in the convergence rates that differ significantly from standard cases. Consequently, the predicted asymptotic behavior becomes observable only for very large computational budgets, which may be impractical. In such scenarios, we propose a mitigation strategy to reconcile theoretical expectations with practical performance.

\section{Control variate method for stochastic differential equations}
\label{sec:main}
We consider a general SDE on the interval $[0,1]$ and study the $\mathbb{R}^d$-valued stochastic process $X_t$ solution of 
\begin{equation}
\label{eq:basic_sde}
\dd X_t = b(X_t) \dd t + \sigma(X_t) \dd W_t, 
\quad  X_0 = \xi,  
\end{equation}
where $W_t$ is a standard Brownian motion, $\xi \in \mathbb{R}^d$, $b ,\sigma :\mathbb{R}^d \to \mathbb{R}^d$ are smooth vector fields (with bounded derivatives).
The objective is to estimate a quantity of interest (QoI) of the form $I=\EE[F(X_1)]$, where $F$ is a smooth and bounded function.

We will first describe a standard MC approach to the estimation of the QoI based on the use of a primary discretization scheme.
We will next describe an original CV approach based on the appropriate use of two discretization schemes (primary and secondary) that have specific properties.

\subsection{Primary discretization scheme}
As we will see below, 
the primary discretization scheme should have the following properties so as to allow for an efficient MC estimation method of the QoI $I$:\\
- it is based on the small time step $h'=1/N'$.\\
- for each time step, it requires a limited number of calls to the drift $b$ (ideally, only one call).\\
- its implementation depends on a family of i.i.d. (independent and identically distributed) standard normal variables ${\bf G}^{{\rm p},h'}$ that is obtained by a linear transform ${\cal G}^{{\rm p},h'} : {\cal C}([0,1],\mathbb{R}) \to \mathbb{R}^{d^{{\rm p},h'}}$ applied to the Brownian motion $(W_t)_{t \in [0,1]}$.
Here the upper index notation ${\rm p}$ stands for primary.
\\
- it should have a high weak order $\alpha$.

In this paper, the total number of calls to the drift $b$ is interpreted as the cost of the method because it is often this term that is computationally expensive (as in the case of Brownian particles in complex potentials).

The primary scheme gives an approximation $X_1^{h'}$ of the solution of the SDE (\ref{eq:basic_sde}) at time $1$ that is a deterministic function of ${\bf G}^{{\rm p},h'}$.
This random vector is made of $d^{{\rm p},h'}$  i.i.d. standard normal variables.
We present below two schemes that can be used as primary schemes.\\

{\bf Euler--Maruyama scheme.}
The solution of (\ref{eq:basic_sde}) can be approximated with the Euler--Maruyama scheme as follows:  $X_0^{h'} = \xi$, for each $1 \leq i' \leq N'$,  $t_{i'} = i'h'$ and 
\begin{equation}
\label{eq:tradi_euler}
X_{t_{i'}}^{h'} 
= X_{t_{i'-1}}^{h'} 
+  b \left ( X_{t_{i'-1}}^{h'} \right ) h' 
+ \sigma \left ( X_{t_{i'-1}}^{h'} \right ) \sqrt{h'} g^{{\rm p}}_{i'}.
\end{equation}
We remark that $ X_{1}^{h'} $
(and more generally, $ \{ X_{t_{i'}}^{h'} ,i'=1,\ldots,N'   \}$) is a deterministic function of 
\begin{equation}
\label{eq:defGph1}
{\bf G}^{{\rm p},h'}= \left\{ g^{{\rm p}}_{i'} ,i'=1,\ldots,N' \right\},
\end{equation}
where 
$$
g^{{\rm p}}_{i'}= \frac{1}{\sqrt{h'}} \big( W_{i'h'}-W_{(i'-1)h'} \big)
$$
are indeed i.i.d. standard normal variables.
At each iteration, the drift coefficient $b(\cdot)$ is called only once.
The method has weak order $\alpha=1$ \cite{kloeden92}.\\

{\bf SRA1 scheme.}
When $\sigma$ is constant, we can consider a higher-order method such as a stochastic RK method SRA1, that is:  $X_0^{{\rm rk},h'} = \xi$, for each $1 \leq i' \leq N'$,  $t_{i'} = i'h'$ and 
\begin{equation}
\label{eq:sra1_part1}
X_{t_{i'}}^{{\rm rk},h'} 
= X_{t_{i'-1}}^{{\rm rk},h'} 
+ \left ( \frac{1}{3} b \left ( X_{t_{i'-1}}^{{\rm rk},h'} \right ) +
\frac{2}{3} b \left ( X_{t_{i'-1}}^{{\rm rk},h'} + \frac{3}{4} \theta_{i'-1} \right )
\right ) h' 
+ \sigma \sqrt{h'} g^{{\rm p}}_{i'}  ,
\end{equation}
where 
\begin{equation}
\label{eq:sra1_part2}
\theta_{i'-1} 
= 
b \left ( X_{t_{i'-1}}^{{\rm rk},h'} \right ) h' + 
\sigma \sqrt{h'} \left (g^{{\rm p}}_{i'} + \frac{1}{\sqrt{3}} f^{{\rm p}}_{i'} \right ).
\end{equation}
We remark that $ X_{1}^{h'} $
is  a deterministic function of 
\begin{equation}
\label{eq:defGph2}
{\bf G}^{{\rm p},h'}= \{ g^{{\rm p}}_{i'}, f^{{\rm p}}_{i'} ,i'=1,\ldots,N' \},
\end{equation}
where 
\begin{align}
g^{{\rm p}}_{i'} =& \frac{1}{\sqrt{h'}} \big( W_{i'h'}-W_{(i'-1)h'} \big), \\
f^{{\rm p}}_{i'} =&  \frac{\sqrt{3}}{\sqrt{h'}}\Big[ 
\frac{2}{h'} \int_{(i'-1)h'}^{i'h'} (W_u-W_{(i'-1)h'} ) {\rm d} u  - 
\big(W_{i'h'}-W_{(i'-1)h'}\big)
\Big]
\label{eq:deffpi}
\end{align}
are indeed i.i.d. standard normal variables.
Here, at each iteration, the drift coefficient $b(.)$ is called twice. 
%We use the notation $g_{i'} = g_i^j$ and $f_{i'} = f_i^j$ when $i' = (i-1)p+j$.
See Equation (6.14) and Table 6.3 in \cite{rossler10} for details on SRA1. This is a two-stage stochastic RK scheme.
The method has weak order $\alpha=2$.

\subsection{Optimized standard Monte-Carlo method}
\label{subsec:MCstandard}
In a standard MC method, one uses the primary discretization method 
%(such as the Euler--Maruyama method or a higher-order Runge--Kutta method) 
with fine time size $h'$ for the SDE that provides an approximation $X_1^{h'}$ of $X_1$, and one estimates the QoI $I=\EE[F(X_1)]$ by
\begin{equation}
\label{eq:st_MC}
\hat{I}_{M'}^{h'} = \frac{1}{M'} \sum_{k'=1}^{M'} F(X_1^{h',k'})
,
\end{equation}
where $X_1^{h',k'}$, $k'=1,\ldots,M'$, are i.i.d. copies of $X_1^{h'}$.

The quadratic error ${\cal E}  = \EE [ (\hat{I}_{M'}^{h'} - I)^2]$ is the sum of a squared bias term and a variance term:
$$
{\cal E} = \big( \EE[F(X_1^{h'})] - \EE[F(X_1)]\big)^2 + \frac{1}{M'} {\rm Var} \big( F(X_1^{h'}) \big).
$$
If the discretization method has weak order $\alpha$, that is to say:
$$
\big| \EE[ F(X_1^{h'}) ] - \EE[F(X_1)]  \big| \lesssim {h'}^\alpha,
$$
then the quadratic error satisfies:
$$
{\cal E}  \lesssim {h'}^{2\alpha} + \frac{1}{M'}  ,
$$
because  ${\rm Var}\big(  F( {X}_1^{h'}) \big)  \simeq {\rm Var}\big(  F({X}_1) \big) (1+o(1))$.
Here $\lesssim$ / $\simeq$ means that the inequality /  equality is valid up to a multiplicative constant.

The total cost of the simulation is 
$$
{\cal C}  \simeq \frac{M'}{h'} ,
$$
because the grid step is $h'$ (hence there are $1/h'$ time steps) and the MC sample size is $M'$.

For a fixed total cost, one can minimize (the upper bound of) the error by choosing 
\begin{equation}
\label{eq:st_optimized_parameters}
h' \simeq {\cal C}^{-\frac{1}{2\alpha+1}} 
\: \mbox{and} \: M' \simeq {\cal C}^{\frac{2\alpha}{2\alpha+1}}
\end{equation}
and then 
\begin{equation}
\label{eq:st_optimized_error}
{\cal E}\lesssim {\cal C}^{-\frac{2\alpha}{2\alpha+1}}.
\end{equation}
This result is obtained by looking at the critical points of the function $h'\mapsto {h'}^{2\alpha}+{\cal C}^{-1} {h'}^{-1}$ obtained by substituting $M'={\cal C} h'$ into the expression of the upper bound of the error ${h'}^{2\alpha} + \frac{1}{M'} $.
This result has been known for a long time \cite{duffie95} and numerical examples there and in \cite{kloeden92} show that this asymptotic behavior can be clearly observed on finite samples.

\begin{remark}
\label{rem:euler}
For the Euler--Maruyama method, one has $\alpha=1$ \cite{kloeden92}  and then ${\cal E}\lesssim{\cal C}^{-\frac{2}{3}}$.
For second-order schemes such as the ones proposed in \cite{milshtein79,talay84,talay90} or the stochastic RK scheme  SRA1 (which is valid when the diffusion coefficient is constant and which requires two calls to the drift function at each step) \cite{rossler10}, one has $\alpha=2$ and then ${\cal E}\lesssim{\cal C}^{-\frac{4}{5}}$.
\end{remark}

%Below there are two versions of the parabolic BM driven ODE. 
%The first one, independent of $\mathcal{G}$, is obtained by sampling directly the coefficients $W_{s,t}$ and $H_{s,t}$ as shown in Equation \eqref{eq:para_bm} with $s = (i-1)h$ and $t = ih$.
%It is used to estimate the expectation of the control variate.
%The second version is (practical) coupled with the Euler scheme \eqref{eq:tradi_euler}, in the sense that they have a common source of randomness $\mathcal{G}$. The random source $\parab{\mathcal{G}} \subset \mathcal{G}$, which appears in  the PPBM but not in the Euler scheme, is introduced to correct the bias from using the composite trapezium rule only. 
%On each interval $[(i-1)h,ih]$ where $1 \leq i \leq N$, the practical parabolic approximation of BM is given by the formula \eqref{eq:pract_para_bm}.
%It is used to build the control variate.

\subsection{Secondary discretization scheme}
As we will see below, the secondary discretization scheme should have the following properties so as to allow for an efficient CV estimation method of the QoI $I$:\\
- it is based on the moderate time step $h =1/N$, which is such that $h=qh'$, with $q$ an integer that is designed to be large.\\
- for each time step, it requires a limited number of calls to the drift $b$ (ideally, one call).\\
- its implementation depends on a family ${\bf G}^{{\rm s},h}$  of i.i.d. standard normal variables that is obtained by a linear transform ${\cal G}^{{\rm s},h} : {\cal C}([0,1],\mathbb{R}) \to\mathbb{R}^{d^{{\rm s},h}}$ applied to the Brownian motion $(W_t)_{t \in [0,1]}$. Here the upper index notation ${\rm s}$ stands for secondary.\\
- it should have a high strong order $\gamma$.

%The first scheme gives an approximation obtained by sampling $d^{{\rm p},h}$  independent and identically distributed standard normal variables that form the Gaussian vector ${\bf G}^{{\rm p},h}$. This approximation denoted $X_1^{h'}$ is, therefore, a deterministic function of ${\bf G}^{{\rm p},h}$.

As the primary and secondary schemes depend on Gaussian vectors extracted linearly from the driving Brownian motion, it is possible to develop two versions of the secondary scheme.\\
- The unconditioned version is the usual scheme obtained by sampling $d^{{\rm s},h}$ i.i.d. standard normal variables that form the Gaussian vector ${\bf G}^{{\rm s},h}$. This gives the approximation $\check{X}_1^{h}$, which is a deterministic function of ${\bf G}^{{\rm s},h}$.\\
- The conditioned version (conditioned to ${\bf G}^{{\rm p},h'}$) is the usual scheme using a sample that has the distribution of 
${\bf G}^{{\rm s},h}$ given ${\bf G}^{{\rm p},h'}$.
This is possible and easy to implement by the Gaussian conditioning theorem as the joint vector $({\bf G}^{{\rm p},h'}, {\bf G}^{{\rm s},h})$ is Gaussian.
The distribution of ${\bf G}^{{\rm s},h}$ given ${\bf G}^{{\rm p},h'}$
has the form ${\bf A} {\bf G}^{{\rm p},h'} + \hat{\bf B} \hat{\bf G}^{{\rm s},h}$, where ${\bf A}$ is a deterministic $d^{{\rm s},h} \times d^{{\rm p},h'}$ matrix, ${\bf B}$ is a deterministic $d^{{\rm s},h} \times \hat{d}^{{\rm s},h}$ matrix, and $\hat{\bf G}^{{\rm s},h}$ is an independent family of $\hat{d}^{{\rm s},h}$ i.i.d. standard normal variables. This gives the approximation $\hat{X}_1^{h}$, which is a deterministic function of ${\bf G}^{{\rm p},h}$ and $\hat{\bf G}^{{\rm s},h}$.
Therefore, ${X}_1^{h'}$ and $\hat{X}_1^{h}$ are correlated via ${\bf G}^{{\rm p},h'}$. This correlation is critical in the CV method proposed in the next section.
We present below two schemes that can be used as secondary scheme.\\

{\bf Milstein method.}
Milstein method is a discretization scheme that has strong  order $\gamma=1$ \cite{milshtein79}.
In the unconditioned version,
$\check{X}_0^{{\rm mi},h} = \xi$ and 
\[
\check{X}_{t_{i}}^{{\rm mi},h} 
= \check{X}_{t_{i-1}}^{{\rm mi},h} 
+ b \left ( \check{X}_{t_{i-1}}^{{\rm mi},h} \right ) h +
\frac{1}{2}  (\sigma \cdot \nabla )\sigma \left ( \check{X}_{t_{i-1}}^{{\rm mi},h} \right ) h \big( (g^{{\rm s}}_i)^2-1 \big)
+ \sigma \left ( \check{X}_{t_{i-1}}^{{\rm mi},h} \right ) \sqrt{h} g^{{\rm s}}_{i}  ,
\]
for each $1 \leq i \leq N$, with $t_{i} = ih$.
We remark that $ \check{X}_{1}^{h} $
is a deterministic function of ${\bf G}^{{\rm s},h}= \left\{ g^{{\rm s}}_{i}, i=1,\ldots,N \right\}$,
where 
$$
g^{{\rm s}}_{i} = \frac{1}{\sqrt{h}} \big( W_{ih}-W_{(i-1)h} \big)
$$
are indeed i.i.d. standard normal variables.\\
In the conditioned version (conditioned to ${\bf G}^{{\rm p},h'}=\{ g_{i'}^{\rm p}, i'=1,\ldots,N'\}$ defined by (\ref{eq:defGph1}) or to ${\bf G}^{{\rm p},h'}= \{ g^{{\rm p}}_{i'}, f^{{\rm p}}_{i'} ,i'=1,\ldots,N' \}$ defined by (\ref{eq:defGph2})), we have $\hat{X}_0^{{\rm mi},h} = \xi$ and 
\[
\hat{X}_{t_{i}}^{{\rm mi},h} 
= \hat{X}_{t_{i-1}}^{{\rm mi},h} 
+ b \left ( \hat{X}_{t_{i-1}}^{{\rm mi},h} \right ) h +
\frac{1}{2}  (\sigma \cdot \nabla )\sigma \left ( \hat{X}_{t_{i-1}}^{{\rm mi},h} \right ) h \big( (g_i^{{\rm s}} )^2-1 \big)
+ \sigma \left ( \hat{X}_{t_{i-1}}^{{\rm mi},h} \right ) \sqrt{h} {g}_{i}^{{\rm s}}  ,
\]
 for each $1 \leq i \leq N$,  
where $$
g_i^{{\rm s}} =  \frac{1}{\sqrt{q}} \sum_{j=1}^q g_{(i-1)q+j}^{{\rm p}} . $$
Note that this method requires to call the gradient of $\sigma$.
In the case where $\nabla \sigma$ is not known analytically, it is possible to use a slightly different scheme in which the term $(\sigma \cdot \nabla )\sigma \left ( \check{X}_{t_{i-1}}^{{\rm mi},h} \right ) h \big( (g^{{\rm s}}_i)^2-1 \big)$ is replaced by $\sigma \left ( \check{X}_{t_{i-1}}^{{\rm mi},h} + h \big( (g^{{\rm s}}_i)^2-1 \big) \sigma \big ( \check{X}_{t_{i-1}}^{{\rm mi},h} \big )\right) -    \sigma \left ( \check{X}_{t_{i-1}}^{{\rm mi},h} \right ) $.\\

{\bf Parabolic approximation.}
As introduced in \cite{lyons20}, the (piecewise) parabolic approximation of a Brownian motion is the (piecewise) second order polynomial that satisfies some particular interpolation relations detailed in Appendix \ref{app:parab1}.
Let us consider a time step $h=1/N$. The piecewise parabolic approximation $\parab{W}^{{\rm s},h}$ of the Brownian motion $W$ on $[0,1]$ is, for $u \in [(i-1)h,ih)$:
\begin{equation}
\parab{W}_u^{{\rm s},h} =
\sum_{j=1}^{i-1} g_j^{{\rm s}} 
+ \frac{u-(i-1)h}{\sqrt{h}} g_i^{{\rm s}} 
+ \frac{\sqrt{3} (u-(i-1)h)(ih-u)}{h^{3/2}} g_i^{{\rm s},\prime} ,
\label{eq:defpiecewiseparab}
\end{equation}
where 
\begin{equation}
\label{eq:gis_unconditional}
g_i^{{\rm s}} = \frac{W_{(i-1)h,ih}}{\sqrt{h}},\qquad
g_i^{{\rm s},\prime}  = \frac{\sqrt{12} H_{(i-1)h,ih}}{\sqrt{h}} ,
\end{equation}
and $W_{s,t}$ and $H_{s,t}$ are defined by
\begin{equation}
\label{eq:defH}
W_{s,t} = W_t - W_s
\: 
\mbox{ and }
\:
H_{s,t} = -\frac{1}{2} W_{s,t} + \frac{1}{t-s} \int_s^t W_{s,u} \dd u.
\end{equation}
Note that the piecewise parabolic approximation of the Brownian motion  depends only on ${\bf G}^{{\rm s},h}=\{
g_i^{{\rm s}},g_i^{{\rm s},\prime} ,i=1,\ldots,N\}$,
which is a family of i.i.d. standard normal variables.

The unconditioned secondary discretization  scheme consists in solving the SDE by replacing the Brownian motion $W$ by its  piecewise parabolic approximation $\parab{W}^{{\rm s},h}$ and by integrating the corresponding ordinary differential equation. 
If it is possible to integrate the ordinary differential equation driven by the piecewise parabolic approximation exactly,  then one gets an approximation denoted $\check{X}^{h}_1$ that has strong convergence order $\gamma=1$ as shown in \cite{lyons20}.
If it is not possible to integrate the ordinary differential equation driven by the piecewise parabolic approximation exactly, then one can use a RK type solver that requires one call to the drift function (and four calls to the diffusivity function) per time step and that gives an approximation denoted $\check{X}_1^{h}$ that has the same strong order $\gamma=1$. This solver is detailed in Appendix~\ref{app:solvesde}.

The conditioned secondary discretization scheme consists in solving the SDE by replacing the Brownian motion $W$ by its conditional piecewise parabolic approximation $\parab{W}^{{\rm s},h}$ (conditioned to ${\bf G}^{{\rm p},h'}$) and by integrating the corresponding ordinary differential equation. The form of the conditional piecewise parabolic approximation depends on the form of the conditioning ${\bf G}^{{\rm p},h'}$:\\
{\it i)} The conditional piecewise parabolic approximation $\parab{W}^{{\rm s},h}$ conditioned to ${\bf G}^{{\rm p},h'}=\{ g_{i'}^{\rm p}, i'=1,\ldots,N'\}$  defined by (\ref{eq:defGph1})
has the form (\ref{eq:defpiecewiseparab}), but 
$g_i^{{\rm s}},g_i^{{\rm s},\prime}$ are now given by (see Appendices \ref{app:A21} and \ref{app:parab3})
\begin{equation}
\label{eq:gis_conditional}
g_i^{{\rm s}} = \frac{1}{\sqrt{q}} \sum_{j=1}^q g_{(i-1)q+j}^{{\rm p}}, \qquad 
g_i^{{\rm s},\prime}=
\frac{\sqrt{3}}{\sqrt{q}}
\left [ 
\sum \limits_{j=1}^{q} \left (1+\frac{1 - 2j}{q} \right ) g_{(i-1)q+j}^{{\rm p}}
+\frac{1}{\sqrt{3q}} 
\hat{g}_i^{{\rm s}} \right ] ,
\end{equation}
and $\hat{\bf G}^{{\rm s},h} =  \{ \hat{g}_i^{{\rm s}} , i=1,\ldots,N\}$ is a family of i.i.d. standard normal variables independent of ${\bf G}^{{\rm p},h'}$ (they are defined by
$\hat{g}_i^{{\rm s}} = \frac{1}{\sqrt{q}} \sum_{j=1}^q f^{{\rm p}}_{(i-1)q+j}$ with $f^{{\rm p}}_{(i-1)q+j}$ given in terms of $W$ by (\ref{eq:deffpi})).\\
%This gives the approximation denoted $\hat{X}^{h}_1$.\\
{\it ii)} The conditional piecewise parabolic approximation $\parab{W}^{{\rm s},h}$ conditioned to ${\bf G}^{{\rm p},h'}=\{ g_{i'}^{\rm p},f_{i'}^{\rm p}, i'=1,\ldots,N'\}$ defined by (\ref{eq:defGph2})
has the form (\ref{eq:defpiecewiseparab}), but
$g_i^{{\rm s}},g_i^{{\rm s},\prime}$ are now completely determined by ${\bf G}^{{\rm p},h'}$ and given by  (see Appendices \ref{app:A22} and \ref{app:parab4})
\begin{equation}
g_i^{{\rm s}} = \frac{1}{\sqrt{q}} \sum_{j=1}^q g_{(i-1)q+j}^{{\rm p}}, \qquad 
g_i^{{\rm s},\prime}=
\frac{\sqrt{3}}{\sqrt{q}}
\left [ 
\sum \limits_{j=1}^{q} \left (1+\frac{1 - 2j}{q} \right ) g_{(i-1)q+j}^{{\rm p}}
+\frac{1}{\sqrt{3}q} 
\sum_{j=1}^q f^{{\rm p}}_{(i-1)q+j}  \right ] .
\end{equation}
In both cases, this gives the approximation denoted $\hat{X}^{h}_1$.

\subsection{Optimized control-variate method}
\label{subsec:optimcv}
CV estimators are a powerful class of reduced variance estimators \cite{glasserman}.
In the proposed CV method for the solution of the SDE (\ref{eq:basic_sde}) one uses \\
1) 
%the primary and secondary discretization methods with fine time step $h'$ and coarse time step $h=ph'$ for the SDE that provides 
an approximation $(\hat{X}_1^{h},X_1^{h'})$ of $X_1$, where $X_1^{h'}$ is obtained with the primary discretization scheme with fine time step $h'$ %(as in the standard MC method) 
and $\hat{X}_1^{h}$ is obtained with the secondary discretization scheme 
%parabolic approximation of the Brownian motion
 conditioned by the Gaussian vector ${\bf G}^{{\rm p},h'}$ used to generate $X_1^{h'}$ (we will take  a coarse time step $h=qh'$ with $q\gg1 $),\\
2)  %a discretization method with coarse time step $h$ for the SDE that provides
 an approximation $\check{X}_1^{h}$ of $X_1$, based on the unconditioned secondary discretization scheme 
% the parabolic approximation of the Brownian motion
(this is a cheap method because $h$ is much larger than $h'$).\\
One estimates the quantity of interest $I=\EE[F(X_1)]$ by
\begin{equation}
\hat{I}_{M,M'}^{h,h'} = \frac{1}{M} \sum_{k=1}^{M} F(\check{X}_1^{h,k})
+
\frac{1}{M'} \sum_{k'=1}^{M'} \big[ F(X_1^{h',k'}) - F(\hat{X}_1^{h,k'})  \big] ,
\label{eq:defIMMp}
\end{equation}
where $\check{X}_1^{h,k}$, $k=1,\ldots,M$, are i.i.d. copies of $\check{X}_1^{h}$, and $(X_1^{h',k'},\hat{X}_1^{h,k'})$, $k'=1,\ldots,M'$, are i.i.d. copies of $(X_1^h,\hat{X}_1^{h'})$.

%\medskip

%\begin{remark}
%Note that there are two versions of the parabolic BM driven ODE in the control-variate method. 
%The first version ($\hat{X}^h$) is coupled with the Euler scheme \eqref{eq:tradi_euler}, in the sense that they have a common source of randomness $\mathcal{G}$. The random source $\parab{\mathcal{G}} \subset \mathcal{G}$, which appears in 
%the PPBM but not in the Euler scheme, is introduced to correct the bias from using the composite trapezium rule only. 
%On each interval $[(i-1)h,ih]$ where $1 \leq i \leq N$, the practical parabolic approximation of BM is given by the formula \eqref{eq:pract_para_bm}.
%It is used to build the control variate (the second term in (\ref{eq:defIMMp})).
%The second version ($\check{X}^h$) %, independent of $\mathcal{G}$,
%is obtained by sampling directly the coefficients $W_{s,t}$ and $H_{s,t}$ as shown in Equation \eqref{eq:para_bm} with $s = (i-1)h$ and $t = ih$.
%%It is used to estimate the expectation of the control variate (the first term in (\ref{eq:defIMMp})).
%\end{remark}

\begin{remark}
We can take the optimal form of the CV estimator:
$$
\hat{J}_{M,M'}^{h,h'} = \frac{\hat{\lambda}_{M'}^{h,h'}}{M} \sum_{k=1}^M F(\check{X}_1^{h,k})
+
\frac{1}{M'} \sum_{k'=1}^{M'} \big[ F(X_1^{h',k'})-  \hat{\lambda}_{M'}^{h,h'} F(\hat{X}_1^{h,k'})  \big],
$$
with
$$
\hat{\lambda}_{M'}^{h,h'} = \frac{
\frac{1}{M'} \sum_{k'=1}^{M'}  F( {X}_1^{h',k'}) F(\hat{X}_1^{h,k'})
-
\Big[ \frac{1}{M'}\sum_{k'=1}^{M'}  F({X}_1^{h',k'}) 
\Big] \Big[\frac{1}{M'}\sum_{k'=1}^{M'} F(\hat{X}_1^{h,k'})\Big]
}
{
\frac{1}{M'} \sum_{k'=1}^{M'}  F(\hat{X}_1^{h,k'})^2
-
\Big[ \frac{1}{M'}\sum_{k'=1}^{M'}  F(\hat{X}_1^{h,k'}) \Big]^2} .
$$
We can gain some advantage (variance reduction) in using this version, but the optimal control coefficient $\lambda$ is very close to one in our framework because the variables $(X_1^h,\hat{X}_1^{h'})$ are strongly correlated as soon as $h' \ll h \ll 1$. As a consequence, there is no significant loss in using $\lambda=1$ instead of $\hat{\lambda}_{M'}$. 
\end{remark}

%\medskip

{\bf Bias.}
Since $\hat{X}_1^h$ and $\check{X}_1^h$ have the same distribution, the bias of the estimator (\ref{eq:defIMMp}) is 
$$
\EE[\hat{I}_{M,M'}^{h,h'}]-I=\EE[F(X_1^{h'})]-\EE[F(X_1)].
$$

%\medskip

{\bf Quadratic error.}
The quadratic error ${\cal E}  = \EE [ (\hat{I}_{M,M'}^{h,h'} - I)^2]$ is the sum of a squared bias term and two variance terms:
$$
{\cal E} = \Big( \EE[F(X_1^{h'})]-\EE[F(X_1)]\Big)^2
+ \frac{1}{M}  {\rm Var}\Big(  F(\check{X}_1^{h}) \Big) 
+
\frac{1}{M'} {\rm Var}\Big( F(X_1^{h'}) - F(\hat{X}_1^{h})
\Big).
$$
If the discretization method for $X_1^{h'}$ has weak order $\alpha$ and strong order $\beta$, that is to say:
$$
\big| \EE[ F(X_1^{h'}) ] - \EE[F(X_1)]  \big|  
\lesssim {h'}^\alpha, 
\qquad 
 \EE\big[ |X_1^{h'} - X_1|^2\big]^{1/2}
 \lesssim {h'}^\beta, 
$$
and if the discretization method for $ \check{X}_1^{h}$ has strong  order $\gamma$, that is to say:
$$
 \EE \big[ | \check{X}_1^{h} - {X}_1|^2 \big]^{1/2}
 \lesssim {h}^{\gamma},
$$
then the quadratic error satisfies
$$
{\cal E}  \lesssim {h'}^{2\alpha} + \frac{1}{M} +\frac{h^{2\gamma} +{h'}^{2\beta}}{M'}   ,
$$
because ${\rm Var}\big(  F(\check{X}_1^{h}) \big)  \simeq {\rm Var}\big(  F({X}_1) \big) (1+o(1))$.
%and { \color{red}${\rm Var}\big( F(X_1^{h'}) - F(\hat{X}_1^{h}) \big) \lesssim \EE[ |\hat{X}_1^{h}-X_1|^2] +\EE[ |X_1^{h'}-X_1|^2] $.}
%{ \color{red}
%\begin{align*}
%{\rm Var}\big( F(X_1^{h'}) - F(\hat{X}_1^{h})\big)
%=  & {\rm Var}\big( F(X_1^{h'}) - F(X_1) \big) +
%{\rm Var}\big( F(X_1^h) - F(X_1) \big)\\
%& - 2 {\rm Cov}\big( F(X_1^{h'}) - F(X_1), \: F( \hat{X}_1^h) - F(X_1)\big).
%\end{align*}
%}
The total cost of the simulation is 
$$
{\cal C}  \simeq \frac{M}{h} +\frac{M'}{h'} ,
$$
because the time steps are $h$ and $h'$ and the sizes of the MC samples are $M$ and $M'$, respectively, in the two terms of Eq.~(\ref{eq:defIMMp}).

\begin{proposition}
\label{prop:error}%
For a fixed total cost, provided $\alpha,\beta >\gamma/(2\gamma+1)$,  the error  is minimized  by choosing 
\begin{equation}
\label{eq:cv_optimized_parameters}
h \simeq  {\cal C}^{-\frac{1}{4\alpha\gamma + 2\alpha+1}}, \: \:
h' \simeq {\cal C}^{-\frac{2\gamma+1}{4\alpha\gamma + 2\alpha+1}}, \: \:
M \simeq {\cal C}^{\frac{4\alpha\gamma + 2\alpha}{4\alpha\gamma + 2\alpha+1}}, \: \:
\mbox{ and } M'\simeq  {\cal C}^{\frac{4\alpha\gamma+2\alpha-2\gamma}{4\alpha\gamma + 2\alpha+1}}
\end{equation} 
(note that $h' \ll h \ll 1$ and $1\ll M'\ll M$), and then:
\begin{equation}
\label{eq:cv_optimized_error}
{\cal E}\lesssim {\cal C}^{-\frac{4\alpha\gamma + 2\alpha}{4\alpha\gamma  +2\alpha+1}}.
\end{equation}
\end{proposition}
The proof of the proposition is given in Appendix~\ref{app:lag1}.
For example, for the Euler--Maruyama method for $X_1^{h'}$, one has $\alpha=1$, $\beta=1/2$ 
%(when $\sigma$ is constant, we actually have $\beta=1$) 
\cite{kloeden92},  and for the parabola method for $\hat{X}_1^{h}$ one has $\gamma=1$ \cite{lyons20}, and then ${\cal E}\lesssim {\cal C}^{-\frac{6}{7}}$.
If we replace the Euler--Maruyama method by the SRA1 scheme (valid when the diffusion coefficient is constant) \cite{rossler10}, one has $\alpha=2$, $\beta=3/2$ and then ${\cal E}\lesssim{\cal C}^{-\frac{12}{13}}$.

\medskip

% 
%\begin{remark}
%In the parabola method one assumes that one knows how to integrate exactly the ordinary differential equation which comes from the original SDE in which the Brownian motion is replaced by a parabola (or, more exactly, a sequence of parabolas defined on each time increment with time step $h$). We then get exactly $\hat{X}_1^h$. 
%%When the diffusion coefficient of the SDE is constant, 
%It is, however, possible to replace this ideal solver by a Runge--Kutta type solver that requires one call to the drift function of the SDE and that gives an approximation $\hat{X}_1^h$ that has the same strong order $\gamma=1$ (cf Appendix~\ref{app:solvesde}).
%\end{remark}
%

%\medskip

\begin{remark}
The parabola method (or the Milstein method) has weak order $\gamma=1$ \cite{lyons20}. Hence, if ones uses this method alone in the standard MC method described in Subsection \ref{subsec:MCstandard}, then one gets a quadratic error that satisfies ${\cal E} \lesssim {\cal C}^{-\frac{2\gamma}{2\gamma+1}}$, that is to say ${\cal E} \lesssim {\cal C}^{-\frac{2}{3}}$.
Comparing with Remark~\ref{rem:euler} (that deals with optimized standard MC methods) and Eq.~(\ref{eq:cv_optimized_error}), it is clear that the proposed CV method has a clear advantage over the optimal uses of the two distinct methods.
\end{remark}

%\medskip
 
%\begin{remark}
%As explained in \cite[pages 341-342]{kloeden92}, the Euler---Maruyama scheme has strong convergence order $0.5$, but when $\sigma$ is constant, it has strong convergence order $1$.
%As shown in  \cite{lyons20}, the parabola method has strong order $1$, even when $\sigma$ is not constant. If, however, $\sigma$ is constant and $\EE[ b''(X_t)]=0$, then the strong order is larger than $1$ (at least, $3/2$).
%\end{remark}

\begin{remark}
The CV approach proposed in this section is compatible with the many-query context.
In this context, the objective is to compute a large number of MC estimations of 
expectations of functions of the solutions of parameterized SDEs \cite{boyaval10}.
This situation is encountered in finance, for instance, where option pricing often involves calculating the price of an option for different values of the parameters of a parameterized local volatility model to be calibrated.
\end{remark}

\section{Numerical illustrations}
\label{sec:num}%
We propose two discretization methods, one based on the Euler method and the other on the RK method as illustrative examples, for the pair $(\hat X^h, X^{h'})$ and $\check X^h$. We recall that $\hat X^h$ is driven by the parabolic approximation of Brownian motion, conditioned on the finer Brownian increments used for $X^{h'}$. The discretization of $\check X^h$ is the same as that of $\hat X^h$, but without conditioning. We recall that the process \((X_t)_{t \geq 0}\) is the solution to Equation~(5). As we will refer to this formulation later on, we present the SDE in its Stratonovich form, that is 
\begin{equation*}
\label{eq:basic_sde_strato}
\dd X_t = \tilde b(X_t) \dd t + \sigma(X_t) \circ \dd W_t, 
\:  X_0 = \xi,  
\end{equation*}
where $\tilde b = b - (1/2)  \sigma' \sigma$. In component notation, for each \( i = 1, \dots, d \), 
\[
\tilde b_i
=
b_i - 
\frac{1}{2} \sum \limits_{j=1}^d \frac{\partial \sigma_i}{\partial x_j} \sigma_j.
\]  

\subsection{Euler scheme based method}
Let \( X_{t_{i'}}^{h'} \) denote the Euler scheme introduced in Equation~\eqref{eq:tradi_euler}, constructed using the fine discretization step \( h' \) and driven by the Brownian increments \( h'^{1/2} g_{i'}^{\rm p} \). We can now formulate the coarse time-step scheme for $\hat{X}^h$. Define the sequence $(\hat{X}_{t_i}^h)_{0 \leq i \leq N}$, with $t_i = i h$, by setting $\hat{X}_0^h = \xi$, and for $1 \leq i \leq N$,  
$
\hat{X}_{t_i}^h = z \left (1; \hat{X}_{t_{i-1}}^h, \hat{A}_{i-1,i} , \hat{B}_{i-1,i} \right )
$
where 
$$
\hat A_{i-1,i} = g_i^{ \rm s}+\sqrt{3} g_i^{{\rm s},\prime},
\: \:
\hat{B}_{i-1,i} = - \sqrt{12} g_i^{{\rm s},\prime}
$$ 
with 
$g_i^{ \rm s}$ and $g_i^{ \rm s, \prime}$  
defined in the  equation \eqref{eq:gis_conditional}
and
for any $z_0,A,B$ and $\{ z(u;z_0,A,B), \: u \in [0,1] \}$ is obtained from
%\vspace{1cm}
\begin{equation*}
\begin{cases}
& z(0) =  z_0,\\
& {\rm d } z(u)/{\rm d } u = h \tilde b(z(u)) + h^{1/2} \sigma(z(u)) (A+Bu) ,\quad  u \in (0,1].
\end{cases}
\end{equation*}
If $z(.)$ does not have a closed form expression, then we propose the approximation with local error of order $O(h^2)$
as follows (see Appendix \ref{app:solvesde} for the derivation of the approximation formula):
$$
z_h(1;z_0,A,B) = z_0+ h B_1+ S_2 - S_0(1-h^{1/2} I_1) + \frac{h^{1/2}}{6} ( S_3 - 2S_1+S_0) I_1,
$$
where we can compute all terms with four calls to $\sigma$ and only one call to $b$: one first computes
$$
S_0 = \sigma(z_0), 
$$
then 
$$
B_1 = b\big(z_0+h^{1/2} S_0 I_3\big) ,\quad S_1 = \sigma\big(z_0 + h^{1/2} S_0 I_1\big) ,
$$
and 
$$
S_2 = \sigma\big( z_0+h S_0 I_2 + h^{3/2} B_1 I_4 \big), \quad S_3= \sigma\big(z_0+h^{1/2} S_0 I_1 +h^{1/2} S_1 I_1 \big).
$$
%{\it Remark.} 
%When $\sigma$ is constant, we only need to call $b$ once: $z_h(1;z_0,A,B)=z_0+h b(z_0+h^{1/2} \sigma I_3)+h^{1/2} \sigma I_1$.
%Here $I_1 = A + B/2$, $I_3 = A/2 + B/6$, and $I_4 = A/2 + B/3$.
%The derivation of the approximation formula is shown in Appendix \ref{app:solvesde}. 
Finally, the coarse time-step scheme for $\check{X}^h$ is defined in a similar way to what was done for $\hat{X}^h$, except that $g_i^{ \rm s}$ and $g_i^{ \rm s, \prime}$  are replaced by an unconditioned version of \eqref{eq:gis_unconditional}.

\subsubsection{Example 1 of an SDE with multiplicative noise}
\begin{figure}[h!]
\centering
\begin{tikzpicture}[scale=0.9]
\begin{loglogaxis}[legend style={at={(1,1)},anchor=north east}, compat=1.3,
  xmin=exp(6.5), xmax=exp(18),ymin=exp(-13),ymax=exp(-1),
  xlabel= {$\mathcal{C}$},
  ylabel= {$\mathcal{E}$}]
\addplot[color=black,thick,mark=square,mark size=1pt] table [x index=0, y index=1]{BSIGMAGENERAL1/13.txt};   
\addplot[color=red,thick,mark=square,mark size=1pt] table [x index=0, y index=1]{BSIGMAGENERAL1/1737.txt}; 
\addlegendentry{standard}
\addlegendentry{control variate}
\end{loglogaxis}
\end{tikzpicture}
\begin{tikzpicture}[scale=0.9]
\begin{loglogaxis}[legend style={at={(0,0)},anchor=south west}, compat=1.3,
  xmin=exp(6.5), xmax=exp(18),ymin=exp(-13),ymax=exp(-1),
  xlabel= {$\mathcal{C}$},
  ylabel= {$\mathcal{E}$}]
  
\addplot[color=blue,thick,mark=square,mark size=1pt] table [x index=0, y index=1]{BSIGMAGENERAL1/1727.txt}; 
\addplot[color=red,thick,mark=square,mark size=1pt] table [x index=0, y index=1]{BSIGMAGENERAL1/1737.txt}; 
\addplot[color=black,thick,mark=square,mark size=1pt] table [x index=0, y index=1]{BSIGMAGENERAL1/1747.txt}; 
\addplot[color=red,thick,mark=o,mark size=1pt] table [x index=0, y index=1]{BSIGMAGENERAL1/2737.txt}; 
\addplot[color=black,thick,mark=o,mark size=1pt] table [x index=0, y index=1]{BSIGMAGENERAL1/2747.txt}; 
\addplot[color=black,thick,mark=triangle,mark size=1pt] table [x index=0, y index=1]{BSIGMAGENERAL1/3747.txt}; 
\addlegendentry{-0.83}
\addlegendentry{\underline{-0.86}}
\addlegendentry{-0.83}
\addlegendentry{-0.72}
\addlegendentry{-0.72}
\addlegendentry{-0.57}
\end{loglogaxis}
\end{tikzpicture}
\caption{Empirical quadratic error $\mathcal{E}$ \eqref{eq:empirical_quadratic_error} ($M'' = 1000$) versus cost $\mathcal{C}$ for the MC estimation of a QoI associated with the SDE \eqref{eq:sde1}. The fine discretisation is based on the Euler-Maruyama method ($\alpha=1$).
The parabola ODE method ($\gamma=1$) is used for the CV method.
\textbf{Left:} 
Comparison between the optimized standard MC expressed in \eqref{eq:st_MC} - \eqref{eq:st_optimized_parameters} - \eqref{eq:st_optimized_error} and the optimized CV MC expressed in \eqref{eq:defIMMp} - \eqref{eq:cv_optimized_parameters} - \eqref{eq:cv_optimized_error}. 
The steepest slope $-6/7$ corresponds to the CV strategy and the other $-2/3$ to the standard optimized MC estimation.
\textbf{Right:}
Comparison between different time step sizes and numbers of MC samples in the CV strategy where $h \sim \mathcal{C}^{-x}$, $h' \sim \mathcal{C}^{-y}$, $M \sim \mathcal{C}^{1-x}$, and $M' \sim \mathcal{C}^{1-y}$. We examine six configurations. 
S1) $x = \frac{1}{7}$, $y = \frac{2}{7}$ (blue squares); \underline{S2}) $x = \frac{1}{7}$, $y = \frac{3}{7}$ (red squares); S3) $x = \frac{1}{7}$, $y = \frac{4}{7}$ (black squares); C1) $x = \frac{2}{7}$, $y = \frac{3}{7}$ (red circles); C2) $x = \frac{2}{7}$, $y = \frac{4}{7}$ (black circles); and T1) $x = \frac{3}{7}$, $y = \frac{4}{7}$ (black triangles). The corresponding slopes are shown in the legend.
These are numerical values obtained by log regression of the empirical quadratic errors.
The steepest slope is observed for $x = \frac{1}{7}$, $y = \frac{3}{7}$, which is consistent with the theory.
}
\label{fig:BSIGMAGENERAL1}
\end{figure}
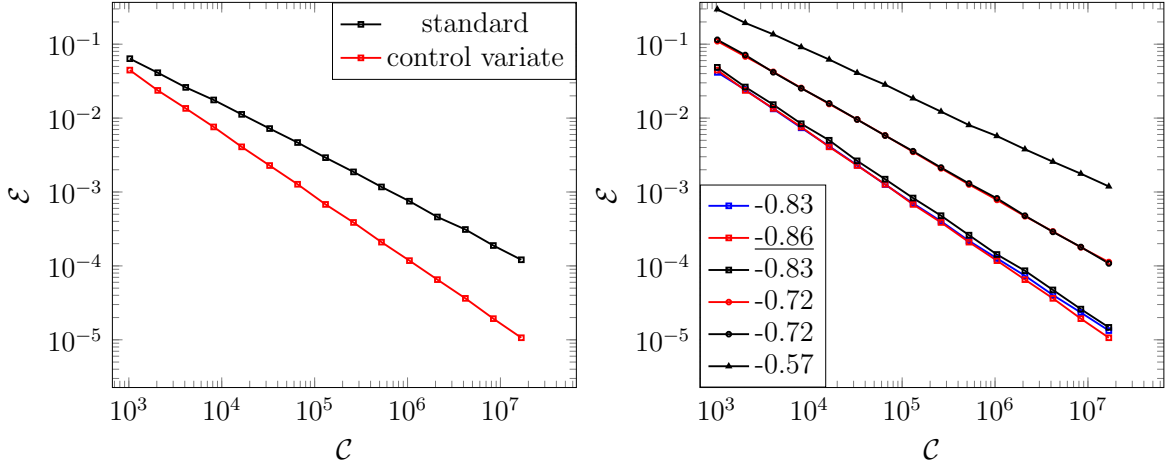

We apply our Euler-based method to an example with non-degenerate multiplicative noise:
\begin{equation}
\label{eq:sde1}
\dd X_t = \frac{1}{2} X_t \dd t + \sqrt{1+ X_t^2} \dd W_t, \, t >0, \, \mbox{with an initial condition} \, X_0 \in \mathbb{R}.
\end{equation}
The equation has an explicit solution
 $$
X_t = \sinh\left(\sinh^{-1}(X_0) + W_t\right) = X_0 \cosh(W_t) + \sqrt{X_0^2 + 1} \sinh(W_t).
 $$
 In the Stratonovich formulation, it takes a driftless form that is easy to integrate.
 Therefore at $t=1$, one has $\EE X_1 = X_0 e^{1/2}$.
The empirical quadratic error associated with the CV estimation of $\EE X_1$ is 
\begin{equation}
\label{eq:empirical_quadratic_error}
\mathcal{E}_{M,M',M''}^{h,h'}
= 
%\EE ( I_{M,M'}^{h,h'} - \EE X_1)^2 
\frac{1}{M''} \sum \limits_{m''=1}^{M''} \left ( \hat{I}_{M,M',m''}^{h,h'} - \EE X_1 \right )^2
\end{equation}
where $\left \{ \hat{I}_{M,M',m''}^{h,h'}, \: 1 \leq m'' \leq M'' \right \}$ is a family of $M''$  independent copies of $\hat{I}_{M,M'}^{h,h'}$, see \eqref{eq:defIMMp}. 
Figure \ref{fig:BSIGMAGENERAL1} illustrates the empirical error ($M''=10^3$) in \eqref{eq:empirical_quadratic_error} w.r.t. the computational cost $\C$ for the SDE \eqref{eq:sde1}. The QoI has a closed-form expression. The left part confirms the theoretical convergence rates for both the standard and CV methods. The right part compares several parameter configurations,  validating the optimal choices $M\sim \mathcal{C}^{6/7}$,
$M'\sim \mathcal{C}^{4/7}$,
$h \sim \mathcal{C}^{-1/7}$ and $h' \sim \mathcal{C}^{-3/7}$ , which lead to the minimized quadratic error $\mathcal{E} \sim \C^{-6/7}$.

\subsubsection{Example 2 of an SDE with multiplicative noise} 
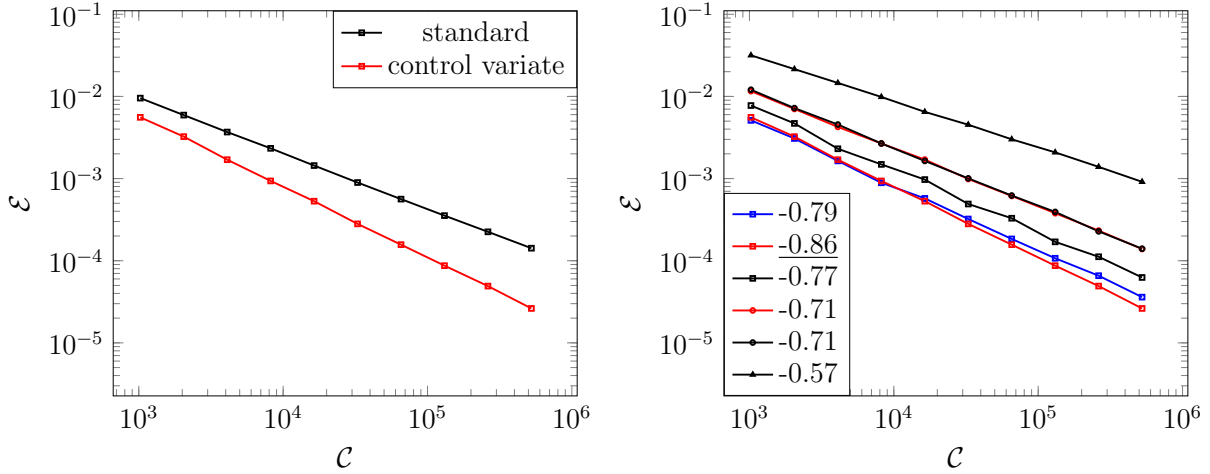
\begin{figure}[h!]
\centering
\begin{tikzpicture}[scale=0.9]
\begin{loglogaxis}[legend style={at={(1,1)},anchor=north east}, compat=1.3,
  %xmin=exp(6.5), xmax=exp(13.9),ymin=exp(-13)+0.00002,ymax=0.0061+0.00002,
  xmin=exp(6.5), xmax=exp(13.9),ymin=exp(-13),ymax=exp(-1)*0.3,
  xlabel= {$\mathcal{C}$},
  ylabel= {$\mathcal{E}$}]
  \addplot[color=black,thick,mark=square,mark size=1pt] table [x index=0, y index=1]{BSIGMAGENERAL2/13.txt}; 
\addplot[color=red,thick,mark=square,mark size=1pt] table [x index=0, y index=1]{BSIGMAGENERAL2/1737.txt}; 
\addlegendentry{standard}
\addlegendentry{control variate}
\end{loglogaxis}
\end{tikzpicture}
\begin{tikzpicture}[scale=0.9]
\begin{loglogaxis}[legend style={at={(0,0)},anchor=south west}, compat=1.3,
  %xmin=exp(6.5), %xmax=exp(13.9),ymin=exp(-13)+0.00002,ymax=0.0061+0.00002,
  xmin=exp(6.5), xmax=exp(13.9),ymin=exp(-13),ymax=exp(-1)*0.3,
  xlabel= {$\mathcal{C}$},
  ylabel= {$\mathcal{E}$}]
  
\addplot[color=blue,thick,mark=square,mark size=1pt] table [x index=0, y index=1]{BSIGMAGENERAL2/1727.txt}; 
\addplot[color=red,thick,mark=square,mark size=1pt] table [x index=0, y index=1]{BSIGMAGENERAL2/1737.txt}; 
\addplot[color=black,thick,mark=square,mark size=1pt] table [x index=0, y index=1]{BSIGMAGENERAL2/1747.txt}; 
\addplot[color=red,thick,mark=o,mark size=1pt] table [x index=0, y index=1]{BSIGMAGENERAL2/2737.txt}; 
\addplot[color=black,thick,mark=o,mark size=1pt] table [x index=0, y index=1]{BSIGMAGENERAL2/2747.txt}; 
\addplot[color=black,thick,mark=triangle,mark size=1pt] table [x index=0, y index=1]{BSIGMAGENERAL2/3747.txt}; 
\addlegendentry{-0.79}
\addlegendentry{\underline{-0.86}}
\addlegendentry{-0.77}
\addlegendentry{-0.71}
\addlegendentry{-0.71}
\addlegendentry{-0.57}
\end{loglogaxis}
\end{tikzpicture}
\caption{Empirical quadratic error $\mathcal{E}$ versus cost $\mathcal{C}$ for the MC estimation of a QoI related to the SDE \eqref{eq:sde2}. 
Same configuration as Figure~\ref{fig:BSIGMAGENERAL1}. The steepest slope is again observed for $x=1/7$, $y=3/7$, which is consistent with the theory.
}
\label{fig:BSIGMAGENERAL2}
\end{figure}

We consider a simple example with non-degenerate multiplicative noise:
\begin{equation}
\label{eq:sde2}
\dd X_t = - X_t \dd t + \sqrt{1+X_t^2} \dd W_t, \, t >0, \, \mbox{with an initial condition} \, X_0 \in \mathbb{R}.
\end{equation}
The solution of this equation does not have an explicit form. However, by the well-known Feynman-Kac formula, the solution to the partial differential equation
\[
u_t(x,t) - x u_x(x,t) + \frac{1 + x^2}{2} u_{xx}(x,t) = 0, \quad \text{for } t < 1,
\]
with terminal condition \( u(x,1) = x \), satisfies
$u(x,t) = \mathbb{E} \left( X_1 \mid X_t = x \right).$
We solve this equation using a standard finite difference method with very high accuracy to approximate \( u(X_0,0) \approx \mathbb{E} X_1 \). This approximation is then used to compute the empirical quadratic error, following the same expression as in~\eqref{eq:empirical_quadratic_error}, except that \( \mathbb{E}X_1 \) is replaced by \( u(X_0,0) \).
Figure \ref{fig:BSIGMAGENERAL2}  is analogous to Figure \ref{fig:BSIGMAGENERAL1} and supports the same conclusions, but applies to the SDE \eqref{eq:sde2}, which does not have a closed-form solution. The QoI is computed using a PDE method. The observed convergence slopes again support the theoretical predictions, further reinforcing the robustness of the CV method through another example involving an SDE with multiplicative noise.

\subsection{RK scheme based method: Brownian particle in a double well potential}
\begin{figure}[h!]
\centering
\begin{tikzpicture}[scale=0.9]
\begin{loglogaxis}[legend style={at={(1,1)},anchor=north east}, compat=1.3,
  xmin=exp(6.5), xmax=exp(13.9),ymin=exp(-13),ymax=exp(-1)*0.3,
  xlabel= {$\mathcal{C}$},
  ylabel= {$\mathcal{E}$}]

\addplot[color=black,thick,mark=square,mark size=1pt] table [x index=0, y index=1]{RKBROWNIANPARTICLE/15.txt}; 
\addplot[color=red,thick,mark=square,mark size=1pt] table [x index=0, y index=1]{RKBROWNIANPARTICLE/113313.txt}; 
\addlegendentry{standard}
\addlegendentry{control variate}
\end{loglogaxis}
\end{tikzpicture}
\begin{tikzpicture}[scale=0.9]
\begin{loglogaxis}[legend style={at={(0,0)},anchor=south west}, compat=1.3,
  xmin=exp(6.5), xmax=exp(13.9),ymin=exp(-13),ymax=exp(-1)*0.3,
  xlabel= {$\mathcal{C}$},
  ylabel= {$\mathcal{E}$}]
 
\addplot[color=blue,thick,mark=square,mark size=1pt] table [x index=0, y index=1]{RKBROWNIANPARTICLE/113213.txt}; 
\addplot[color=red,thick,mark=square,mark size=1pt] table [x index=0, y index=1]{RKBROWNIANPARTICLE/113313.txt}; 
\addplot[color=black,thick,mark=square,mark size=1pt] table [x index=0, y index=1]{RKBROWNIANPARTICLE/113413.txt}; 

\addplot[color=red,thick,mark=o,mark size=1pt] table [x index=0, y index=1]{RKBROWNIANPARTICLE/213313.txt}; 
\addplot[color=black,thick,mark=o,mark size=1pt] table [x index=0, y index=1]{RKBROWNIANPARTICLE/213413.txt}; 

\addplot[color=black,thick,mark=triangle,mark size=1pt] table [x index=0, y index=1]{RKBROWNIANPARTICLE/313413.txt}; 

\addlegendentry{-0.77}
\addlegendentry{\underline{-0.93}}
\addlegendentry{-0.90}
\addlegendentry{-0.85}
\addlegendentry{-0.85}
\addlegendentry{-0.77}
\end{loglogaxis}
\end{tikzpicture}
\caption{Empirical quadratic error $\mathcal{E}$ \eqref{eq:empirical_quadratic_error} versus cost $\mathcal{C}$ for the MC estimation of a QoI associated with the SDE \eqref{eq:sde3}. The fine discretization is based on the stochastic RK method (SRA1) ($\alpha=2$).
The parabola ODE method ($\gamma=1$) is used for the CV method.
\textbf{Left:} 
Comparison between the optimized standard MC expressed in \eqref{eq:st_MC} - \eqref{eq:st_optimized_parameters} - \eqref{eq:st_optimized_error} and the optimized CV MC expressed in \eqref{eq:defIMMp} - \eqref{eq:cv_optimized_parameters} - \eqref{eq:cv_optimized_error}. 
The steepest slope $-12/13$ corresponds to the CV strategy and the other $-4/5$ to the standard optimized MC estimation.
\textbf{Right:}
Comparison between different time step sizes and number of MC samples in the CV strategy where $h \sim \mathcal{C}^{-x}$, $h' \sim \mathcal{C}^{-y}$, $M \sim \mathcal{C}^{1-x}$, and $M' \sim \mathcal{C}^{1-y}$. We examine six configurations. 
S1) $x = \frac{1}{13}$, $y = \frac{2}{13}$ (blue squares); \underline{S2}) $x = \frac{1}{13}$, $y = \frac{3}{13}$ (red squares); S3) $x = \frac{1}{13}$, $y = \frac{4}{13}$ (black squares); C1) $x = \frac{2}{13}$, $y = \frac{3}{13}$ (red circles); C2) $x = \frac{2}{13}$, $y = \frac{4}{13}$ (black circles); and T1) $x = \frac{3}{13}$, $y = \frac{4}{13}$ (black triangles). The corresponding slopes are shown in the legend.
The steepest slope is observed for $x = \frac{1}{13}$, $y = \frac{3}{13}$, which is consistent with the theory.
}
\label{fig:RKBROWNIANPARTICLE3}
\end{figure}
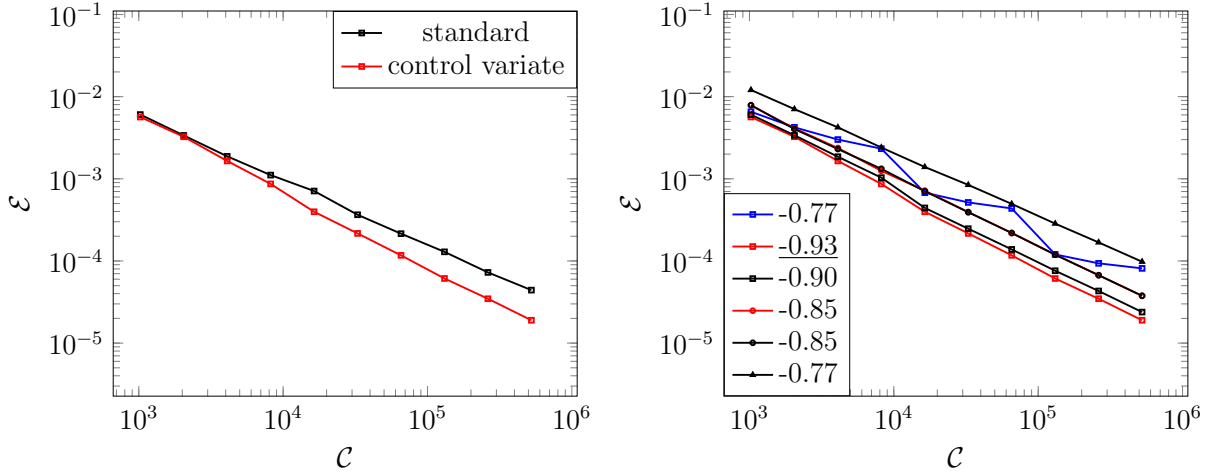

We consider a Brownian particle in a double well potential whose derivative is
$
U'(x) = x(x+1)(x-2). 
$
Thus, we define $b(x) = - U'(x)$ and the Brownian particle has the following dynamics
\begin{equation}
\label{eq:sde3}
\dd X_t = b(X_t) \dd t + \sigma \dd W_t, \, t >0, \, \mbox{with an initial condition} \, X_0 \in \mathbb{R}.
\end{equation}
There is no explicit formula for $(X_t)_{t \geq 0}$. 
Figure \ref{fig:RKBROWNIANPARTICLE3} presents results for the double-well potential SDE \eqref{eq:sde3} using the SRA1 scheme \eqref{eq:sra1_part1}-\eqref{eq:sra1_part2} as the primary method. The results are structured in a manner similar to the previous two figures. The SDE does not have a closed-form solution and thus the QoI is computed using a PDE method. In the left part, the steepest slope of $-{12}/{13}$ (i.e., $\mathcal{C}^{-12/13}$) demonstrates the theoretical advantage of the CV method over the standard MC approach, which is characterized by a gentler slope of $-{4}/{5}$ (i.e., $\mathcal{C}^{-4/5}$). The right part compares several parameter configurations,  validating the optimal choices $M\sim \mathcal{C}^{12/13}$,
$M'\sim \mathcal{C}^{10/13}$,
$h \sim \mathcal{C}^{-1/13}$ and $h' \sim \mathcal{C}^{-3/13}$ , which lead to the minimized quadratic error $\mathcal{E} \sim \C^{-12/13}$.

\subsection{Special case: Inhomogeneous geometric Brownian motion}
We consider the inhomogeneous geometric Brownian motion (IGBM) that is the solution of 
\begin{equation*}
\label{eq:igbm}
\dd X_t = a (b-X_t) \dd t + \sigma X_t \dd W_t, \, t >0, \, \mbox{with an initial condition} \, X_0 \in \mathbb{R}.
\end{equation*}
This process has no known method of exact simulation, it is frequently encountered in mathematical finance and other fields, and it is used as a model to test and compare various numerical methods \cite{tubikanec}.
We recall there is an explicit formula, using the notation $\tilde a = a + \sigma^2/2$, that is
\begin{equation}
\label{eq:igbm_expression}
\forall t \geq 0, \: X_t = 
e^{- \tilde a t + \sigma W_t} \left ( y_0 + ab \int_0^t e^{ \left (\tilde a u - \sigma W_u \right )} \dd u \right ).
\end{equation}
Therefore, with $N \geq 1, \, h := 1/N$ fixed as above, for any $1 \leq i \leq N$,
$$
X_{ih} = e^{- \tilde a h + \sigma (W_{ih} - W_{(i-1)h}) } \left ( X_{(i-1)h} + ab \int_{(i-1)h}^{ih} e^{\tilde a (u-(i-1)h) - \sigma (W_u-W_{(i-1)h})} \dd u \right ).
$$
% We will use the parabola ODE method as proposed in \cite{lyons20} for the IGBM.
% It boils down to replacing $w$ by its parabolic approximation on the interval $[(i-1)h,ih]$ (in the Stratonovich form).
We will build our CV estimator using the following discretization schemes. 
\begin{enumerate}
\item \textbf{Euler--Maruyama method with fine time step $h' = 1/N'$.}\\
Define the sequence $(X_{t_{i'}'}^{h'})_{0\leq i' \leq N'}$ where $t_{i'}' := i'h'$ by $X_0^{h'} := X_0$ and
$$
X_{t_{i'}'}^{h'} =
X_{t_{i'-1}'}^{h'}
+ a (b-X_{t_{i'-1}'}^{h'}) h'
+ \sigma X_{t_{i'-1}'}^{h'} \sqrt{h'} g_{i'}^{\rm p}, \: \mbox{for} \: 1 \leq i' \leq N'.
$$

%Here $(i,j)$ is the unique pair such that $k'= (i-1)q + j$, $1 \leq i \leq N$, and $1 \leq j \leq q$. 

%We recall that the $g_{i}^j$ are the normalized elementary increments in $\mathcal{G}_0$.

\item \textbf{Conditioned-parabola ODE method with coarse time step $h$.}\\
Define the sequence $(\hat{X}_{t_i}^h)_{0\leq i \leq N}$ where $t_i := i h$ by $\hat{X}_0^h := X_0$ and
$$
\hat{X}_{t_i}^h = e^{- \tilde a h + \sigma \sqrt{h} g_{i}^{\rm s} } \left ( \hat{X}_{t_{i-1}}^h + ab \int_{(i-1)h}^{ih} e^{\tilde a (u-(i-1)h) - \sigma 
\left ( \parab{W}_u^{{\rm s},h} - \parab{W}_{(i-1)h}^{{\rm s},h}  \right )} \right ) \dd u ,
$$
where $\parab{W}_u^{{\rm s},h}$
is defined in \eqref{eq:defpiecewiseparab}, depending on $g_i^{\rm s}$ and $g_i^{\rm s,\prime}$ in  
\eqref{eq:gis_conditional}.
We remark
$$
\int_{(i-1)h}^{ih} e^{\tilde a (u-(i-1)h) - \sigma \left ( \parab{W}_u^{{\rm s},h} - \parab{W}_{(i-1)h}^{{\rm s},h}  \right )} \dd u 
= \mathcal{I}(h,\hat{W}_{i-1,i}, \hat{H}_{i-1,i}),
$$
where 
$$
\mathcal{I}(h,W,H) := h \int_0^1 e^{\tilde a h v - \sigma (v W + 6 v (1-v)H)} \dd v.
$$
and 
$$
\hat{W}_{i-1,i} = \sqrt{h} g_i^{\rm s}
\mbox{ and } 
\hat{H}_{i-1,i}
= \frac{\sqrt{h}}{2\sqrt{3}} g_i^{\rm s, \prime}.
$$
\item \textbf{Free-parabola ODE method with coarse time step $h$.}\\
Define the sequence $(\check{X}_{t_i}^h)_{0\leq i \leq N}$ where $t_i := i h$ by $\check{X}_0^h := X_0$ and
$$
\check{X}_{t_i}^h = e^{- \tilde a h + \sigma \check{W}_{i-1,i} } \left ( \check{X}_{t_{i-1}}^h + ab \: \mathcal{I}(h,\check{W}_{i-1,i},\check{H}_{i-1,i}) \right ),
$$
where 
$$
\check{W}_{i-1,i} = \sqrt{h} g_i^{ \rm s} 
\mbox{ and }
\check{H}_{i-1,i} = \frac{\sqrt{h}}{2 \sqrt{3}} g_i^{ \rm s, \prime}.
$$
Here, the variables $g_i^{ \rm s}$ and $g_i^{ \rm s, \prime}$ are an unconditioned version of \eqref{eq:gis_unconditional}.
\end{enumerate}
The integral $\mathcal{I}(h,W,H)$ can be approximated by a 3-point Gauss-Legendre quadrature numerical method as follows
$$
\mathcal{I}(h,W,H)
\approx
\frac{h}{18} 
\left \{ 5 g \left (  \frac{1}{2} -\frac{1}{2} \sqrt{\frac{3}{5}}   \right ) 
+ 8 g \left ( \frac{1}{2}  \right ) 
+ 5 g \left ( \frac{1}{2} + \frac{1}{2} \sqrt{\frac{3}{5}}  \right ) \right \},
$$
where 
$
g(v) = e^{\tilde a h v - \sigma (v W + 6 v (1-v)H)}
$
depends on $(h,W,H)$.
Using Eq.~\eqref{eq:igbm_expression} at $t=1$, one has $\EE X_1 = \exp(-a) X_0 + b (1-\exp(-a))$.   

\begin{figure}
\centering
\begin{tikzpicture}[scale=0.7]
\begin{axis}[
    xlabel={$\log(h)$},
    ylabel={$E(h)$},
    grid=none,
    legend style={at={(0.35,0.85)},anchor=north,legend columns=-1},
    width=12cm,
    height=8cm
]
\addplot[color=blue, mark=*] coordinates {
(-15.24923797, -20.32745959) (-14.55609079, -19.22882282) (-13.86294361, -18.38156157) (-13.16979643, -17.61948329) (-12.47664925, -16.89354661) (-11.78350207, -16.18429973) (-11.09035489, -15.48310879) (-10.39720771, -14.78623053) (-9.70406053, -14.08951337) (-9.01091335, -13.39480207) (-8.31776617, -12.70205962) (-7.62461899, -12.01104459) (-6.93147181, -11.31528514) (-6.23832463, -10.6193062) (-5.54517744, -9.92059228) (-4.85203026, -9.22255158) (-4.15888308, -8.54250407) (-3.4657359, -7.79867887) (-2.77258872, -7.06756805) (-2.07944154, -6.34096193) (-1.38629436, -5.53185512) (-0.69314718, -4.70604549) (0.0, -3.90605627) 
};
\addlegendentry{$\log(E(h))$}

\addplot[color=red, mark=square*] coordinates {
(-15.24923797, -19.87359931) (-14.55609079, -19.15207207) (-13.86294361, -18.43054482) (-13.16979643, -17.70901758) (-12.47664925, -16.98749034) (-11.78350207, -16.2659631) (-11.09035489, -15.54443586) (-10.39720771, -14.82290862) (-9.70406053, -14.10138138) (-9.01091335, -13.37985414) (-8.31776617, -12.65832689) (-7.62461899, -11.93679965) (-6.93147181, -11.21527241) (-6.23832463, -10.49374517) (-5.54517744, -9.77221793) (-4.85203026, -9.05069069) (-4.15888308, -8.32916345) (-3.4657359, -7.60763621) (-2.77258872, -6.88610896) (-2.07944154, -6.16458172) (-1.38629436, -5.44305448) (-0.69314718, -4.72152724) (0.0, -4.0) 
};
\addlegendentry{$ \log(h) - 4$}

\end{axis}
\end{tikzpicture}
\caption{Strong convergence error $  E(h) = \EE \big[ | \check{X}_1^{h} - {X}_1|^2 \big]^{1/2}$ for the parabola method and the IGBM case. We numerically find that the strong error is  $\eps h^\gamma$ with $\gamma=1$ (as predicted by the general theory) and $\eps=\exp(-4)$.}
\label{fig:IGBM_parabola_strong}
\end{figure}
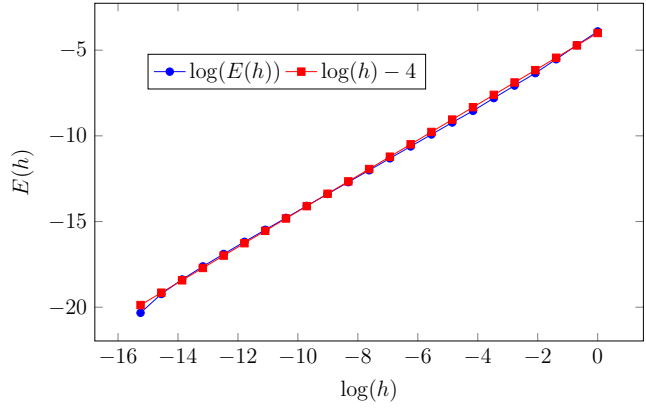

Before presenting the numerical results we discuss in the following remark a special property of the IGBM.

\begin{remark}
\label{rem:6}
The strong error for the parabola method is of order one.
Figure~\ref{fig:IGBM_parabola_strong} shows that the strong error for the IGBM is of the form 
$ \EE \big[ | \check{X}_1^{h} - {X}_1|^2 \big]^{1/2}
 \simeq \eps {h}^{\gamma}$ with $\gamma=1$ and the multiplicative constant $\eps \simeq \exp(-4)$ is approximately one hundred times smaller than the corresponding constant  for the strong error of the Euler--Maruyama scheme. 
In other words, we actually deal with a situation in which the quadratic error satisfies 
\begin{equation}
\label{epsE}
{\cal E}  \lesssim {h'}^{2\alpha} + \frac{1}{M} +\frac{\eps^2 h^{2\gamma} +{h'}^{2\beta}}{M'} .
\end{equation}
In Appendix~\ref{app:lag2}
we revisit the proof of Proposition~\ref{prop:error} to determine the optimal parameters $h,h',M,M'$ by keeping track of the small constant $\eps$. We assume $\beta \geq 1/2$. If ${\cal C} \gg \eps^{-4\alpha}$, then we have
\begin{align*}
&
h \simeq  {\cal C}^{-\frac{1}{4\alpha\gamma + 2\alpha+1}} \eps^{-\frac{4\alpha}{4\alpha\gamma + 2\alpha+1}} , \: \:
h' \simeq {\cal C}^{-\frac{2\gamma+1}{4\alpha\gamma + 2\alpha+1}}
\eps^{\frac{2}{4\alpha\gamma + 2\alpha+1}}, \\
& M \simeq {\cal C}^{\frac{4\alpha\gamma + 2\alpha}{4\alpha\gamma + 2\alpha+1}}\eps^{-\frac{4\alpha}{4\alpha\gamma + 2\alpha+1}}, \: \:
 M'\simeq  {\cal C}^{\frac{4\alpha\gamma+2\alpha-2\gamma}{4\alpha\gamma + 2\alpha+1}}\eps^{\frac{2}{4\alpha\gamma + 2\alpha+1}} ,
\end{align*}
and the error is
$$
{\cal E}\lesssim {\cal C}^{-\frac{4\alpha\gamma + 2\alpha}{4\alpha\gamma  +2\alpha+1}} \eps^{\frac{4\alpha}{4\alpha\gamma + 2\alpha+1}}.
$$
If ${\cal C} \ll \eps^{-4\alpha}$, then we have
$$
h \simeq 1, \: \:
h' \simeq {\cal C}^{-\frac{1}{2\alpha+1}}
\eps^{\frac{2}{2\alpha+1}}, \: \:
M\simeq {\cal C}, \: \:
M' \simeq {\cal C}^{\frac{2\alpha}{2\alpha+1}}
\eps^{\frac{2}{2\alpha+1}},
$$
and the error satisfies
$$
{\cal E}\lesssim {\cal C}^{-1}.
$$
\end{remark}

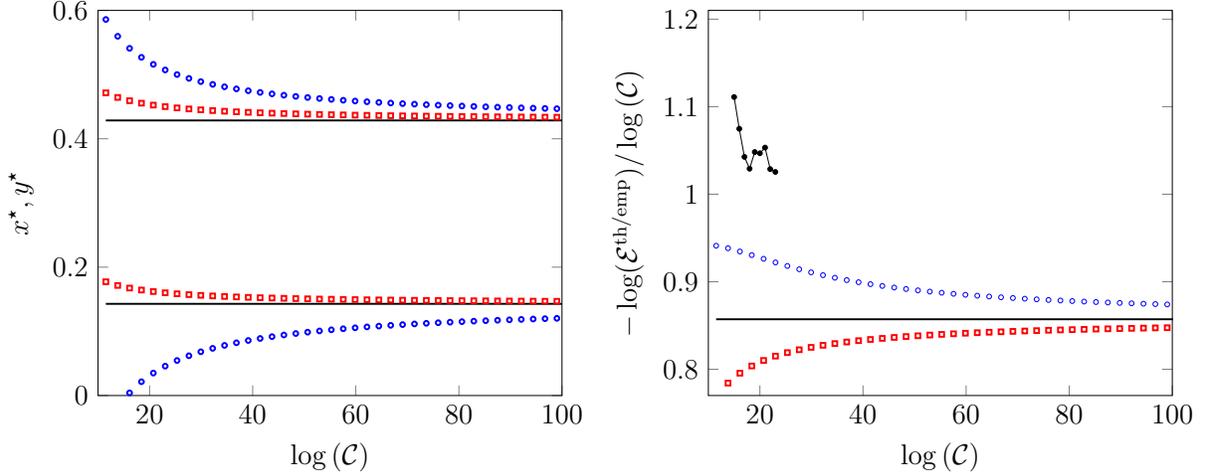
\begin{figure}[ht]
\centering
\begin{tikzpicture}[scale=0.9]
\begin{axis}[
  legend style={at={(0,1)},anchor=north west},
  xmin=10, xmax=100,
  ymin=0, ymax=0.60,
  xlabel= {$\log \left (\mathcal{C} \right )$},
  ylabel= {$x^\star,y^\star$}
]
\addplot[only marks,thick,color=red,mark=square,mark size=1pt] table [x index=0, y index=1] {IGBM/IGBM_theoretical_minimizers.txt};%x epsilon = 1
\addplot[only marks,thick,color=blue,mark=o,mark size=1pt] table [x index=0, y index=2] {IGBM/IGBM_theoretical_minimizers.txt};%x epsilon = 0.01
\addplot[thick,color=black,mark size=1pt] table [x index=0, y index=3] {IGBM/IGBM_theoretical_minimizers.txt}; % limit
\addplot[only marks,thick,color=red,mark=square,mark size=1pt] table [x index=0, y index=4] {IGBM/IGBM_theoretical_minimizers.txt};
\addplot[only marks,thick,color=blue,mark=o,mark size=1pt] table [x index=0, y index=5] {IGBM/IGBM_theoretical_minimizers.txt};
\addplot[thick,color=black,mark size=1pt] table [x index=0, y index=6] {IGBM/IGBM_theoretical_minimizers.txt};
\end{axis}
\end{tikzpicture}
\begin{tikzpicture}[scale=0.9]
\begin{axis}[
  legend style={at={(0,1)},anchor=north west},
  xmin=10, xmax=100,
  ymin=0.77, ymax=1.21,
  xlabel= {$\log \left (\mathcal{C} \right )$},
  ylabel= {$-\log (\mathcal{E}^{\rm th /emp}) /\log \left (\mathcal{C} \right )$}
]
\addplot[only marks, thick,color=red,mark=square,mark size=1pt] table [x index=0, y index=1] {IGBM/IGBM_theoretical_error.txt};
\addplot[only marks, color=blue,mark=o,mark size=1pt] table [x index=0, y index=2] {IGBM/IGBM_theoretical_error.txt};
\addplot[thick,color=black,mark size=0.5pt] table [x index=0, y index=3] {IGBM/IGBM_theoretical_error.txt};
\addplot[color=black,mark=*,mark size=1pt] table [x index=0, y index=1] {IGBM/IGBM_empirical_error.txt};
\end{axis}
\end{tikzpicture}
\caption{
This figure illustrates the comparison between the minimized theoretical errors and the empirical errors observed in simulation.
\textbf{Left:} 
$x^\star$ below $0.2$ and $y^\star$ above $0.4$
versus $\log(\mathcal{C})$. 
These parameters minimize the right-hand-side of \eqref{epsE}, where $h,h',M$ and $M'$ take the form $h \sim \mathcal{C}^{-x^\star}, h' \sim \mathcal{C}^{-y^\star}, M \sim \C^{1-x^\star}$ and $M' \sim \mathcal{C}^{1-y^\star}$.
\textbf{Right:}  $-\log(\mathcal{E}^{\rm th})/\log(\mathcal{C})$ where $\mathcal{E}^{\rm th}$ is the  minimized theoretical quadratic error involving the parameter \( \epsilon \). 
The black dots linked by straight lines correspond to  $ -\log(\mathcal{E}^{\rm emp})/\log(\mathcal{C})$ versus $\log(\mathcal{C})$ where the empirical error $\mathcal{E}^{\rm emp}$ similar to \eqref{eq:empirical_quadratic_error} is estimated with $h \sim \mathcal{C}^{-x^\star}, h' \sim \mathcal{C}^{-y^\star}, M \sim \C^{1-x^\star}$ and $M' \sim \mathcal{C}^{1-y^\star}$. Due to the computational time limitation, we only show $\log (\mathcal{C}) \in [15,20]$. In both subfigures, the range of values for $\log(\mathcal{C})$ is $[10, 100]$. The curves marked with circles corresponds to $\epsilon = 10^{-2}$, while the curves marked with squares corresponds to $\epsilon = 1$. The solid line is the asymptotic limit for $\C$ large.
}
\label{fig:IGBM_-logE/logC}
\end{figure}

Remark~\ref{rem:6} helps to understand the results in Figure \ref{fig:IGBM_-logE/logC} when ${\cal C}$ is not very large.
The left part indicates how the optimal parameters \( x^\star \) and \( y^\star \) evolve with respect to \( \log(\C) \). These parameters minimize the bound on the theoretical quadratic error involving the parameter \( \epsilon \), where $h,h',M$ and $M'$ take the form $h \sim \mathcal{C}^{-x^\star}, h' \sim \mathcal{C}^{-y^\star}, M \sim \C^{1-x^\star}$ and $M' \sim \mathcal{C}^{1-y^\star}$. The bound is shown in Equation \eqref{epsE}. This evolution highlights the transition between pre-asymptotic and asymptotic regimes.
%Indeed one explores the impact of the small constant \( \epsilon \) in the strong convergence rate for the IGBM case. 
In the right part, the graph shows that the empirical error (zigzagging black points) (estimated with $h \sim \mathcal{C}^{-x^\star}, h' \sim \mathcal{C}^{-y^\star}, M \sim \C^{1-x^\star}$ and $M' \sim \mathcal{C}^{1-y^\star}$) we observe from simulations tends to follow qualitatively the same pattern as the minimized theoretical error (curve marked with circles) when $\eps = \exp(-4)$. In both parts, for sufficiently large values of $\C$, the effect of $\eps$ becomes negligible.

\section{Towards a Multilevel Monte Carlo approach}
%and concluding comments}
\label{sec:multi}%
The multilevel Monte Carlo (MLMC) method is a variance reduction technique for simulating SDEs. Instead of relying solely on fine discretizations (which are accurate but expensive), MLMC combines simulations across a hierarchy of time steps: many cheap simulations with coarse discretizations and fewer simulations with fine discretizations. The main idea is to estimate expectations using a telescoping sum of corrections between successive levels. This greatly reduces the overall computational cost while maintaining the accuracy of fine-level simulations \cite{giles08a,giles08b,giles15}.

Let $L \in \mathbb{N}^\star$ be the number of numerical schemes or numerical levels used in MLMC generalizing our estimator \eqref{eq:defIMMp} for which $L=2$. 
 We use the notation $\bh = (h_1,\ldots,h_L)$  where $0 < h_1<\cdots < h_L$ and $\bM = (M_1,\ldots,M_L) \in (\mathbb{N}^\star)^L$.  
Define 
\begin{equation}
\label{eq:defIMMpML}
I_{\bh}^{\bM} =  
\frac{1}{M_L} \sum \limits_{m_L=1}^{M_L} F(\check X_1^{h_L, m_L})
 + \sum \limits_{\ell = 1}^{L-1} \frac{1}{M_\ell} \sum \limits_{m_\ell=1}^{M_\ell}
F(X_1^{h_\ell, m_\ell})
- F(\hat X_1^{h_{\ell+1}, m_\ell}),
\end{equation}
where we consider \( L \) independent families of i.i.d. random variables: 
\[
\left\{ \check{X}^{h_L,m_L} \right\}_{m_L=1}^{M_L}, \quad \text{and for each } \ell = 1, \dots, L-1, \text{ the pairs } 
\left\{ \left( \hat{X}^{h_{\ell+1},m_{\ell}}, X^{h_\ell,m_\ell} \right) \right\}_{m_\ell=1}^{M_\ell}.
\]

In the multilevel estimator \eqref{eq:defIMMpML} one uses \\
1) 
$L-1$ independent approximations $(\hat{X}_1^{h_{\ell+1}},X_1^{h_\ell})$ of $X_1$ for each $1 \leq \ell \leq L-1$ , where $X_1^{h_{\ell}}$ is obtained with a local primary discretization scheme with the time step $h_{\ell}$ 
finer than $h_{\ell+1}$
and $\hat{X}_1^{h_{\ell+1}}$ is obtained with a local secondary discretization scheme 
%parabolic approximation of the Brownian motion
 conditioned by the Gaussian increments used to generate $X_1^{h_\ell}$. 
 %(we will take  a coarse time step $h=qh'$ with $q\gg1 $)
Note that the local secondary scheme employed $\hat{X}_1^{h_{\ell+1}}$ at level $\ell$ which is employed as the local primary scheme for $X_1^{h_{\ell+1}}$ at level $\ell+1$. Consequently, both variables $\hat{X}^{h_{\ell+1}}$ and $X^{h_{\ell+1}}$ are generated using the same numerical scheme and share the same distribution. However, $\hat{X}^{h_{\ell+1}}$ is conditioned on Gaussian vector used to generate $X^{h_\ell}$, whereas $X^{h_{\ell+1}}$ is generated independently (unconditioned).
 \\
2)
 an approximation $\check{X}_1^{h_L}$ of $X_1$, based on the unconditioned local secondary discretization scheme of the last level $L$. 
 
{\bf Bias.}
Since for each \(1 \leq \ell \leq L-1\), the random variables \(\hat{X}_1^{h_\ell}\) and \(X_1^{h_\ell}\), as well as \(\hat{X}_1^{h_L}\) and \(\check{X}_1^{h_L}\), share the same distribution, the bias of the estimator (\ref{eq:defIMMpML}) is (telescoping terms) 
$$
\EE[I_{\bh}^{\bM}]-I=\EE[F(X_1^{h_1})]-\EE[F(X_1)].
$$

{\bf Quadratic error.}
The quadratic error ${\cal E}  = \EE [ (I_{\bh}^{\bM} - I)^2]$ is the sum of a squared bias term and the variance term:
$$
{\cal E} = \Big( \EE[F(X_1^{h_1})]-\EE[F(X_1)]\Big)^2
+ {\rm Var}\big(  I_{\bh}^{\bM} \big). 
$$

The variance term can be decomposed as:
\[
{\rm Var}\big(I_{\bh}^{\bM}\big) = \frac{1}{M_L} {\rm Var}\left( F(\check{X}_1^{h_L}) \right)
+ \sum_{\ell=1}^{L-1} \frac{1}{M_\ell} {\rm Var}\left( F(X_1^{h_\ell}) - F(\hat{X}_1^{h_{\ell+1}}) \right).
\]
If the discretization method for $X_1^{h_1}$ has weak order $\alpha$
and if, for each $1 \leq \ell \leq L$, the discretization method for $ \hat{X}_1^{h_\ell}$ (or $X_1^{h_\ell}$) has strong  order $\gamma_\ell$
then the quadratic error satisfies
\begin{equation}
\label{eq:error_multilevel}
{\cal E}  \lesssim {h_1}^{2\alpha} + \frac{1}{M_L} 
+
\sum \limits_{\ell=1}^{L-1}
\frac{h_\ell^{2\gamma_\ell} +{h_{\ell+1}}^{2\gamma_{\ell+1}}}{M_\ell}   ,
\end{equation}
where we use ${\rm Var}\big(  F(\check{X}_1^{h_L}) \big)  \simeq {\rm Var}\big(  F({X}_1) \big) (1+o(1))$.
The total budget cost is 
$$
\C = \sum \limits_{\ell=1}^L \frac{M_\ell}{h_\ell}. 
$$
In essence, the sum on the right-hand side of the inequality \eqref{eq:error_multilevel} is asymptotically dominated by the term that decays to zero most slowly when $h_L$ and $M_\ell$ go to $0$ and $+\infty$, respectively.
This implies that, in the asymptotic regime, the leading contributions to the error arise from only two levels, effectively reducing the  analysis to the case 
$L=2$.

\appendix

\section{Parabolic approximation of Brownian motion}

This appendix provides the mathematical foundation for the parabolic approximation of Brownian motion, including its interpolation properties and Gaussian conditioning. It presents explicit formulas essential for implementing the numerical schemes discussed in the main text.
\subsection{Definition of the parabolic approximation}
\label{app:parab1}%
We consider a real-valued standard Brownian motion $\{ W_u, u \geq 0 \}$. Fix $t>s$.
%, we work on the interval $[s,t]$. 
As introduced in \cite{lyons20}, the parabolic approximation of $W$ on $[s,t]$, denoted by $\parab{W}$, is the second order polynomial that satisfies 
$$
\parab{W}_s = W_s, \: 
\parab{W}_t = W_t, \: \mbox{and} \:
\int_s^t \parab{W}_u \dd u = \int_s^t W_u \dd u.
$$
It is uniquely defined by the three conditions above, which lead to the following formula
$$
\parab{W}_u 
= W_s 
+ \frac{u-s}{t-s} W_{s,t}
+ \frac{6 (u-s)(t-u)}{(t-s)^2} H_{s,t},
$$
where $W_{s,t}$ and $H_{s,t}$ are defined by (\ref{eq:defH}).
As a consequence, 
$$
\int_s^t \parab{W}_u \dd u =
(t-s) \left ( W_s 
+ \frac{1}{2}W_{s,t}
+ H_{s,t}
\right )
\: \: \mbox{ and } \: \:
\frac{\dd}{\dd u}\parab{W}_u = \frac{1}{t-s} W_{s,t}
+ \frac{6(t+s -2 u)}{(t-s)^2} 
H_{s,t}.
$$
It is clear that $(W_{s,t},H_{s,t})$ is a centered Gaussian vector. A direct calculation shows that $\E ( W_{s,t} H_{s,t} ) = 0$. Therefore, $W_{s,t}$ and $H_{s,t}$  are independent, centered, 
$
\E W_{s,t}^2 = t-s
\: \mbox{ and } \:
\E H_{s,t}^2 = (t-s)/12.
$
We can therefore write
\begin{equation*}
\label{eq:para_bm}
\parab{W}_u
=
W_s 
+ \frac{u-s}{\sqrt{t-s}} g
+ \frac{\sqrt{3} (u-s)(t-u)}{(t-s)^{3/2}} g'.
\end{equation*}
where
\begin{equation*}
g = \frac{W_{s,t}}{\sqrt{t-s}},\qquad
g' = \frac{\sqrt{12} H_{s,t}}{\sqrt{t-s}}
\end{equation*}
are independent and identically distributed standard normal variables.

%Thus
%\begin{align*}
%\frac{\dd \parab{W}_u^{t-s,g,g'}}{\dd u}
%& =
%\frac{g}{(t-s)^{1/2}} + \frac{\sqrt{3}}{(t-s)^{3/2}} (t+s-2u) g'.\\
%%%& =
%%\underbracket{\frac{g}{h^{1/2}} + \frac{\sqrt{3}}{h^{3/2}} 
%%(t+s) g'}_{A =} 
%%+ \underbracket{\left ( - \frac{2 \sqrt{3} g'}{h^{3/2}} %\right )}_{B =} u.
%\end{align*}

\subsection{Two-level time discretization and piecewise parabolic Brownian motion}
In this section we explain the dependence structure of the increments of the Brownian motions and the coefficients of the parabolic approximation. The results follow from Gaussian conditioning theorem.
Consider $N' \in \mathbb{N}^\star$.
% and $\mathcal{G} = \{ g_1, \dots, g_{N'} \}$ an i.i.d. family of $\mathcal{N}(0,1)$.
Let $h' = 1/N'$, $q \in \mathbb{N}^\star$. The coarse time step is $h= p h'$, and the corresponding number of time steps is $N = N'/q$ (provided $N'$ is a multiple of $q$).

\subsubsection{Decomposition of Brownian motion (part 1)}
\label{app:A21}
We denote 
$\delta W_i^j=W_{(i-1)h +jh'}-W_{(i-1)h+(j-1)h'}$, $\delta {\bf W}_i=(\delta W_i^j)_{j=1}^q$, 
$\Delta W_i = W_{ih}-W_{(i-1)h}$ and 
$\Delta I_i = \int_{(i-1)h}^{ih} W_u-W_{(i-1)h} {\rm d} u$.\\
{\it Result.}
The distribution of 
$(\Delta W_i,\Delta I_i)_{i=1}^N$ given $(\delta {\bf W}_i)_{i=1}^N$ is Gaussian with mean 
$$
\Big(
\sum_{j=1}^q \delta W_i^j ,
 h' \sum_{j=1}^q(q+\frac{1}{2}-j)  \delta W_i^j  
\Big)_{i=1}^N
$$
and covariance
$$
\oplus_{i=1}^N 
\begin{pmatrix} 0 & 0\\
0& \frac{qh'^3}{12}  \end{pmatrix}.
$$
In other words,  given $(\delta {\bf W}_i)_{i=1}^N$, 
the vector $(\Delta W_i )_{i=1}^N$ is deterministic and equal to $(\sum_{j=1}^q \delta W_i^j )_{i=1}^N$ and the vector $(\Delta I_i )_{i=1}^N$ is made of independent  Gaussian variables with means $ h' \sum_{j=1}^q(q+\frac{1}{2}-j)  \delta W_i^j  $ and variances $\frac{qh'^3}{12} $.\\
{\it Proof.}
We first note that $(\delta {\bf W}_i,\Delta W_i,\Delta I_i)_{i=1}^N$ are independent and identically $(q+2)$-dimensional Gaussian vectors.
Their distribution is
$$
\begin{pmatrix}
\delta {\bf W}_i\\
\Delta W_i\\
\Delta I_i
\end{pmatrix} 
\sim
{\cal N} \bigg(
\begin{pmatrix}
{\bf 0}\\
0\\
0
\end{pmatrix} ,
\begin{pmatrix}
h' {\bf I}_q & h' {\bf 1} & h'^2 (j+\frac{1}{2}-q)_{j=1}^q \\
h' 
{\bf 1}^T & h & \frac{h^2}{2}\\
 h'^2 {(j+\frac{1}{2}-q)_{j=1}^q}^T & \frac{h^2}{2} & \frac{h^3}{3}
\end{pmatrix} 
\bigg)
$$
We then apply the Gaussian conditioning theorem to get the distribution of $(\Delta W_i,\Delta I_i)$ given $\delta {\bf W}_i$ for any $i$. \qed

We denote by $g_i^j$ the normalized elementary increments $$g_i^j=(W_{(i-1)h+jh'}-W_{(i-1)h+(j-1)h'})/\sqrt{h'}.$$
From the previous result, we find that,
given the elementary increments  $(g_i^j)_{1\leq j \leq q, 1\leq i \leq N}$,
$W_{ih}-W_{(i-1)h}$ is deterministic:
\begin{align*}
W_{ih}-W_{(i-1)h} = \sqrt{h'} G_i , 
\end{align*}
with $G_i=  \sum_{j=1}^q g_i^j$, $W_{ih}$ is deterministic
$$
W_{ih} = \sqrt{h'}\sum_{i'=1}^i G_{i'}
$$
and $\int_{(i-1)h}^{ih} W_u du$  is
of the form 
\begin{align*}
\int_{(i-1)h}^{ih} W_u {\rm d} u 
= 
{h'}^{3/2} \Big( q \sum_{i'=1}^{i-1} G_{i'} +
\sum_{j=1}^q(q+\frac{1}{2}-j) g_i^j +
\frac{\sqrt{q}}{2 \sqrt{3}}  \hat{G}_i \Big),
\end{align*}
where the $\hat{G}_i$'s are independent standard Gaussian variables (and independent from the $g_i^j$'s). 

\subsubsection{Decomposition of Brownian motion (part 2)}
\label{app:A22}
Define
${{\bf H}_i} = (H_i^j)_{j=1}^q$ where
$$
H_i^j = -\frac{1}{2} \delta W_i^j 
+ \frac{1}{h'} \int_{(i-1)h + (j-1)h'}^{(i-1)h + jh'} ( W_u- W_{(i-1)h + (j-1)h'} )\dd u.
$$
The pair 
$(\Delta W_i,\Delta I_i)_{i=1}^N$ given $(\delta {\bf W}_i)_{i=1}^N$ and $( {{\bf H}_i})_{i=1}^N$ has the following deterministic representation
$$
\Delta W_i 
= \sum \limits_{j=1}^q \delta W_i^j
\: 
\mbox{ and } 
\:
\Delta I_i = h' \sum \limits_{j=1}^q \left ( H_i^j + \left ( q-j+\frac{1}{2} \right )\delta W_i^j \right ).
$$
We denote by $f_i^j$ the normalized elementary ``parabolic adjustment" defined by
$$
f_i^j = \frac{2 \sqrt{3}}{{h'}^{1/2}} H_i^j.  
$$
We can notice that the normalized elementary increment $g_i^j = \delta W_i^j/{h'}^{1/2}$ and ``parabolic adjustment" $f_i^j$ are independent, both distributed with a standard normal distribution.

\subsubsection{Conditioned piecewise parabolic Brownian motion (part 1)}
%\label{sec:PPBM}
\label{app:parab3}
We define
$
\mathcal{G}_0 = 
\{ g_{i}^j, \: 
1 \leq i \leq N, \: 
1 \leq j \leq q \}.
$
The lower index runs through the grid with the coarse time step while the upper index is for the grid with the finer one. We use the notation  
$$
G_{i} = \sum \limits_{j=1}^q g_{i}^j.
$$
Additionally, we also define
$
\hat{\mathcal{G}} = 
\{ \hat {G}_{i}, \: 
1 \leq i \leq N \}$ and $\mathcal{G} = \mathcal{G}_0 \cup \hat{\mathcal{G}}$. The two families $\mathcal{G}_0$ and $\hat{\mathcal{G}}$ are independent and they consist of i.i.d. standard normal variables.
For each $1 \leq i \leq N$, 
the parabolic approximation of Brownian motion (conditioned on $\mathcal{G}$) on $[(i-1)h,ih]$, denoted by $\hat{W}$, is the second-order polynomial that satisfies 
$$
\hat{W}_{(i-1)h} = \sqrt{h'} \sum_{r=1}^{i-1} G_r, \: \: 
\hat{W}_{ih} = \sqrt{h'} \sum_{r=1}^i G_r
$$
and
$$
\int_{(i-1)h}^{ih} \hat{W}_u \dd u = 
(h')^{3/2} \left ( q \sum_{r=1}^{i-1} G_r
+ \frac{1}{2}G_i
+ \sum_{j=1}^{q-1} \left ( q - j \right ) g_i^j
+ \frac{1}{2}\sqrt{\frac{q}{3}} \hat{G}_i \right ).
$$
It is uniquely defined by the three conditions above, which  lead to the following formula
%The practical BM becomes for 
\begin{equation}
\label{eq:pract_para_bm}
\forall u \in [(i-1) h, i h], \: 
\hat{W}_u 
= \sqrt{h'} \sum_{r=1}^{i-1} G_r 
+  (u/h-i+1) \hat{W}_{i-1,i}
+  6(u/h-i+1)(i - u/h) \hat{H}_{i-1,i}
\end{equation}
where
\begin{equation}
\label{practicalWH}
\hat{W}_{i-1,i} = \sqrt{h'} G_i
\: 
\mbox{ and } 
\:
\hat{H}_{i-1,i} =
\frac{\sqrt{h'}}{2}
\left [ 
\sum \limits_{j=1}^{q} \left (1+\frac{1 - 2j}{q} \right ) g_{i}^j
+\frac{1}{\sqrt{3q}} 
\hat{G}_i \right ].
\end{equation}
Direct calculations show that the   pair 
$(\hat{W}_{i-1,i},
\hat{H}_{i-1,i})$
and the  pair 
$(W_{s,t}, H_{s,t})$, with $s = (i-1)h$ and $t = ih$, have the same law.

\subsubsection{Conditioned piecewise parabolic Brownian motion (part 2)}
%\label{sec:PPBM++}
\label{app:parab4}
We define 
$
\mathcal{F}_0 = 
\{ f_{i}^j, \: 
1 \leq i \leq N, \: 
1 \leq j \leq p \}
$
and
$\mathcal{H} = \mathcal{G}_0 \cup \mathcal{F}_0$.
For each $1 \leq i \leq N$, 
the parabolic approximation of Brownian motion (conditioned on $\mathcal{H}$) on $[(i-1)h,ih]$, denoted by $\hat{W}_u^C$, is the second-order polynomial that satisfies 
$$
\hat{W}_{(i-1)h}^C = \sqrt{h'} \sum_{r=1}^{i-1} G_r, \: \: 
\hat{W}_{ih}^C = \sqrt{h'} \sum_{r=1}^i G_r
$$
and
$$
\int_{(i-1)h}^{ih} \hat{W}_u^C \dd u = 
(h')^{3/2} \left ( q \sum_{r=1}^{i-1} G_r
+ \frac{1}{2}G_i
+ \sum_{j=1}^{q-1} \left ( q - j \right ) g_i^j
+ \frac{1}{2\sqrt{3}} \sum \limits_{j=1}^q f_i^j \right ).
$$
It is uniquely defined by the three conditions above, which lead to a similar formula as \eqref{eq:pract_para_bm}, except that 
$\hat{H}_{i-1,i}$ is replaced by 
$\hat{H}_{i-1,i}^C$ with
\begin{equation*}
\label{practicaldoublehatH}
\hat H_{i-1,i}^C =
\frac{\sqrt{h'}}{2}
\left [ 
\sum \limits_{j=1}^{q} \left (1+\frac{1 - 2j}{q} \right ) g_{i}^j
+\frac{1}{\sqrt{3} q} 
\sum  \limits_{j=1}^q f_i^j \right ].
\end{equation*}

\section{Error minimization}
\label{app:lag}%
This appendix contains the technical derivation of the error minimization strategy for the CV MC estimator. It includes the proof of Proposition \ref{prop:error} and an  analysis that accounts for small multiplicative constants in pre-asymptotic regimes.

\subsection{Proof of Proposition \ref{prop:error}}
\label{app:lag1}%
The problem is to minimize
${h'}^{2\alpha} + \frac{1}{M} +\frac{h^{2\gamma} +{h'}^{2\beta}}{M'} $ under the constraint $\frac{M}{h} +\frac{M'}{h'}={\cal C}$.
By substituting $M =  {\cal C} h - \frac{M' h}{h'}$ into the expression of the function to be minimized, one obtains the unconstrained minimization problem: find the minimum of ${\cal H}: (h,h',M') \mapsto {h'}^{2\alpha} + \frac{h'}{{\cal C} h h' - M' h} +\frac{h^{2\gamma} +{h'}^{2\beta}}{M'} $.
A critical point satisfies the three equations $\partial_h {\cal H} = \partial_{h'}{\cal H}=\partial_{M'} {\cal H}=0$. Solving this system of three equations gives the three values for $h,h',M'$.

Two remarks:\\
1) The critical point is a minimum that corresponds to an admissible value for $M'\simeq  {\cal C}^{\frac{4\alpha\gamma+2\alpha-2\gamma}{4\alpha\gamma + 2\alpha+1}}$ (a large integer) because $4\alpha \gamma +2\alpha -2\gamma >0$ (indeed, this is equivalent to $\alpha > \gamma/(2\gamma+1)$).\\
2) We observe that the term ${h'}^{2\beta}$ is negligible in front of $h^{2\gamma}$. Indeed the ratio of these two numbers is $h^{4\gamma \beta +2\beta-2\gamma}$ which is small because $4\gamma \beta +2\beta-2\gamma>0$ (indeed, this is equivalent to $\beta > \gamma/(2\gamma+1)$).

\subsection{Error minimization when $\C$ and $\eps$ are fixed}
\label{app:lag2}
We are interested in minimizing 
${h'}^{2\alpha} + \frac{1}{M} +\frac{\eps^2 h^{2\gamma} +{h'}^{2\beta}}{M'} $
with 
$$
h' = \C^{-y}, \: h = \C^{-x}, \: M' = \C^{1-y}, \: M = \C^{1-x}.
$$
Thus we consider the function
\begin{equation}
\label{eq:theo_error}
e(x,y) = \C^{-2\alpha y} + \C^{x-1} 
+ \C^{y-1}(\eps^2 \C^{-2\gamma x} + \C^{-2 \beta y}).
\end{equation}
At a critical point $(x^\star,y^\star)$, we have
$$
\frac{ \partial e_\eps(x^\star,y^\star)}{\partial x} = 0 
\iff
\C^{x^\star-1} = 2 \gamma \epsilon^2 \C^{y^\star-1-2 \gamma x^\star}
$$
and (with $\beta = 1/2$)
$$
\frac{ \partial e_\eps(x^\star,y^\star)}{\partial y} = 0 
\iff
2\alpha \C^{-2\alpha y^\star} = \epsilon^2 \C^{y^\star-1 - 2\gamma x^\star}.
$$
Taking the $\log$ and dividing by $\log (\C)$ gives
\[
\begin{pmatrix}
1 + 2\gamma & -1 \\
-2\gamma & 1 + 2\alpha
\end{pmatrix}
\begin{pmatrix}
x^\star \\
y^\star
\end{pmatrix}
=
\begin{pmatrix}
\displaystyle \frac{\log(2\gamma \epsilon^2)}{\log (\C)} \\
\displaystyle 1 + \frac{\log\left(\frac{2\alpha}{\epsilon^2}\right)}{\log (\C)}
\end{pmatrix}.
\]
%Then
%\begin{align*}
%\begin{pmatrix}
%x^\star\\
%y^\star
%\end{pmatrix}
%=
%\frac{1}{1+2\alpha+4 \alpha \gamma} 
%\begin{pmatrix}
%1 + 2\alpha & 1 \\
%2\gamma & 1 + 2\gamma
%\end{pmatrix}
%\begin{pmatrix}
%\displaystyle \frac{\log(2\gamma \epsilon^2)}{\log (\C)} \\
%\displaystyle 1 + %\frac{\log\left(\frac{2\alpha}%{\epsilon^2}\right)}{\log (\C)}
%\end{pmatrix}.
%\end{align*}
That is 
$$
x^\star
=
\left ( 1+
\frac{\log( (2 \gamma)^{2 \alpha+1} 2 \alpha) }{\log(\C)} + \frac{\log(\eps^{4 \alpha})}{\log(\C)} \right ) / (1+2 \alpha + 4 \alpha \gamma)
$$
and
$$
y^\star 
=
\left ( 1+ 2\gamma + 
\frac{\log( (4 \gamma \alpha)^{2 \gamma} 2 \alpha) }{\log(\C)} - \frac{\log(\eps^2)}{\log(\C)} \right ) / (1+2 \alpha + 4 \alpha \gamma).
$$

\section{Runge-Kutta-type solver for the parabolic scheme}
This appendix introduces a RK-type solver tailored for ordinary differential equations driven by the parabolic approximation of Brownian motion. It demonstrates that the solver achieves strong convergence of order one, requiring only one evaluation of the drift function and four evaluations of the diffusion function per time step.

\label{app:solvesde}
Let $P:[0,1]\to \mathbb{R}$ be a real-valued function (it is a second-order polynomial in the application). Let $b,\sigma:\mathbb{R}^p \to \mathbb{R}^p$ be smooth functions.
In this appendix  $z(u) \in \mathbb{R}^p$ is the solution to 
\begin{equation}
\label{eq:odeapp1}
\dot z(u) = h b(u) + h^{1/2}\sigma(z(u))P(u), \quad  u>0,
\end{equation}
starting from $z(0) = z_0$.

\begin{lemma}
\label{lem:rk}
Denote 
\begin{align*}
I_1 = 
\int_0^1 P(u) \dd u, \qquad 
I_3  = 
\int_0^1 \int_0^u P(u_1) \dd u_1 \dd u.
\end{align*}
We have
\begin{equation}
z(1) = z_0+ h B_1+ S_2 - S_0(1-h^{1/2} I_1) + \frac{h^{1/2}}{6} ( S_3 - 2S_1+S_0) I_1
+O(h^2)  ,
\label{eq:lemrk}
\end{equation}
where 
$$
S_0 = \sigma(z_0), 
$$
$$
B_1 = b\big(z_0+h^{1/2} S_0 I_3\big) ,\quad S_1 = \sigma\big(z_0 + h^{1/2} S_0 I_1\big) ,
$$
$$
S_2 = \sigma\big( z_0+h S_0 I_2 + h^{3/2} B_1 I_4 \big), \quad S_3= \sigma\big(z_0+h^{1/2} S_0 I_1 +h^{1/2} S_1 I_1 \big).
$$
\end{lemma}
This lemma shows that it is possible to compute an approximation of $z(1)$ accurate up to a term of order $O(h^2)$ with four calls to $\sigma$ and only one call to $b$.
When $\sigma$ is constant, Eq.~(\ref{eq:lemrk}) reads: $z(1)=z_0+h b(z_0+h^{1/2} \sigma I_3)+h^{1/2} \sigma I_1 +O(h^2)$.

The strong order of convergence of the parabolic scheme using the solver presented in Lemma~\ref{lem:rk} at each step then follows from \cite[Theorem 11.5.1]{kloeden92}. This result is similar to the one stated in \cite[Theorem 3.17]{lyons20} that addresses the parabolic scheme in the case where Eq.~(\ref{eq:odeapp1}) can be integrated exactly.

{\it Proof.}
If $f:\mathbb{R}^p \to \mathbb{R}^p$ is sufficiently smooth, 
then for any $u \in [0,1]$,
\begin{align*}
f(z(u)) = & f(z_0) + h^{1/2} (f' \sigma)(z_0) \int_0^u P(u_1) \dd u_1\\
& + h \left ( (f' b)(z_0) u +  [(f' \sigma)'\sigma ](z_0) \int_0^u \int_0^{u_1}  P(u_2) P(u_1)  \dd u_2  \dd u_1 \right ) + O(h^{3/2}),
\end{align*}
where $f' = (\partial f_i / \partial z_j)_{ij}$ stands for the Jacobian matrix of $f$ in $\mathbb{R}^{p \times p}$.
Applying this result with $f=b$ and $f=\sigma$, we have  
$$
%\begin{equation}
%\label{eq:csq1}
\int_0^1 b(z(u)) \dd u 
= b(z_0) 
+  h^{1/2} (b' \sigma)(z_0) \int_0^1 \int_0^u P(u_1) \dd u_1 \dd u 
+ O(h)
$$
%\end{equation}
and
\begin{align*}
\nonumber
& \int_0^1 \sigma(z(u)) P(u) \dd u = \sigma(z_0) \int_0^1 P(u) \dd u  + h^{1/2} (\sigma' \sigma)(z_0) \int_0^1\int_0^u P(u_1) \dd u_1 P(u) \dd u
\nonumber\\
& + h \left ( (\sigma' b)(z_0) \int_0^1 u P(u) \dd u +  [(\sigma' \sigma)'\sigma ](z_0) \int_0^1 \int_0^u \int_0^{u_1}  P(u_2) P(u_1)  P(u) \dd u_2  \dd u_1 \dd u\right )
%\nonumber\\ &
+ O(h^{3/2}).
\label{eq:csq2}
\end{align*}
Therefore, we have 
\begin{equation}
\label{eq:formula_zeta}
z(1) = z_0 + h^{1/2}\zeta_{1/2}  + h \zeta_1  + h^{3/2} \zeta_{3/2}  + O(h^2) ,
\end{equation}
where
\begin{align*}
\zeta_{1/2}  & = \sigma(z_0) I_1 , \\
\zeta_1  & = b(z_0) + \sigma'\sigma(z_0) I_2 , \\
\zeta_{3/2}  & = (b'\sigma)(z_0) I_3 +  (\sigma'b)(z_0) I_4 + [(\sigma'\sigma)'\sigma](z_0) I_5 ,
\end{align*}
and
$I_1 = 
\int_0^1 P(u) \dd u$, $I_2  = 
\int_0^1 \int_0^u P(u_1) P(u) \dd u_1 \dd u$,
$I_3  = 
\int_0^1 \int_0^u P(u_1) \dd u_1 \dd u$,
$ I_4  = 
\int_0^1 u P(u) \dd u$,
$I_5  =  \int_0^1 \int_0^u \int_0^{u_1} P(u_2) P(u_1) P(u) \dd u_2 \dd u_1 \dd u$.
Note that all these quantities depend only on $I_1$ and $I_4$ since $I_2=I_1^2/2$, $I_3=I_1-I_4$, $I_5=I_1^3/6$.
Moreover, we have
\begin{align*}
&hb(z_0)+h^{3/2} b'\sigma(z_0) I_3 = hb\big(z_0+h^{1/2} \sigma(z_0) I_3\big) +O(h^2) ,\\
    &h^{1/2} \sigma(z_0) I_1+h \sigma'\sigma(z_0) I_2+h^{3/2} \sigma' b(z_0) I_4 \\
    &\quad = 
\sigma\big( z_0+h \sigma(z_0) I_2 + h^{3/2} b\big(z_0+h^{1/2} \sigma(z_0) I_3\big) I_4 \big) - \sigma(z_0) \big(1-h^{1/2} I_1\big) +O(h^2) ,  
\end{align*}
(here we use $b\big(z_0+h^{1/2} \sigma(z_0) I_3\big)$ instead of $b(z_0)$ in the first term of the right-hand side in order to save one call to $b$),
and
\begin{align*}
h^{3/2} \left( (\sigma'\sigma)' \sigma \right)(z_0) I_1^3
&= h^{1/2} \Big[
\sigma\Big( z_0 + h^{1/2} \sigma(z_0) I_1 + h^{1/2} \sigma\big(z_0 + h^{1/2} \sigma(z_0) I_1\big) I_1 \Big) \\
&\quad - 2 \sigma\Big( z_0 + h^{1/2} \sigma(z_0) I_1 \Big)
+ \sigma(z_0)
\Big] I_1 + O(h^2).
\end{align*}
Substituting into (\ref{eq:formula_zeta})  gives the desired result.
\qed

\end{document}